\documentclass[aap,MSNbibl,nameyear,seceqn,dvips]{arximspdf}

\doi{10.1214/10-AAP727}
\volume{21}
\issue{4}
\pubyear{2011}
\firstpage{1365}
\lastpage{1399}

\makeatletter
\newtheorem{Theorem}{Theorem}[section]
\newproclaim{Definition}{Definition}[section]
\newtheorem{Proposition}{Proposition}[section]

\newtheorem{Lemma}{Lemma}[section]
\newtheorem{Corollary}{Corollary}[section]
\newproclaim{Remark}{Remark}[section]

\def\bptnote#1{}

\makeatother

\begin{document}
\begin{frontmatter}

\title{Optimal multiple stopping time problem}

\runtitle{Multiple stopping}

\begin{aug}
\author[A]{\fnms{Magdalena} \snm{Kobylanski}\ead[label=e1]{magdalena.kobylanski@univ-mlv.fr}},
\author[B]{\fnms{Marie-Claire}~\snm{Quenez}\corref{}\ead[label=e2]{quenez@math.jussieu.fr}}\\
\and \author[C]{\fnms{Elisabeth}~\snm{Rouy-Mironescu}\ead[label=e3]{Elisabeth.Mironescu@ec-lyon.fr}}

\runauthor{M. Kobylanski, M.-C. Quenez and E. Rouy-Mironescu}
\affiliation{Universit\'e Paris Est, Universit\'e Paris-Diderot and
Ecole Centrale de Lyon}
\address[A]{M. Kobylanski\\
CNRS---UMR 8050 (LAMA)\\
Universit\'e de Marne-la-Vall\'ee\\
5, boulevard Descartes\\
Cit\'e Descartes---Champs-sur-Marne\\
77454 Marne-la-Vall\'ee cedex 2\\
France\\
\printead{e1}} 
\address[B]{M.-C. Quenez\\
CNRS---UMR 7599 (LPMA)\\
Universit\'e Denis Diderot (P7)\\
175 rue du Chevaleret\\
75251 Paris Cedex 05\\
France\\
\printead{e2}}
\address[C]{E. Rouy-Mironescu\\
CNRS \\
Ecole Centrale de Lyon\\
Institut Camille Jordan\\
Universit\'e de Lyon\\
43, boulevard du 11 novembre 1918\\
69622 Villeurbanne Cedex\\
France\\
\printead{e3}}
\end{aug}

\received{\smonth{10} \syear{2009}}
\revised{\smonth{6} \syear{2010}}

%
\begin{abstract}
We study the optimal multiple stopping time problem defined for each
stopping time $S$ by $
v(S)=\operatorname{ess}\sup_{\tau_1,\ldots,\tau_d \geq S} E[\psi(
\tau_1,\ldots,\tau_d)|\mathcal{F}_S]$.

The key point is the construction of a \textit{new reward} $\phi$ such
that the value function $v(S)$ also satisfies $
v(S)=\operatorname{ess}\sup_{\theta\geq S } E[\phi(\theta
)|\mathcal{F}_S]$. This
new reward $\phi$ is not a right-continuous adapted process as in the
classical case, but a family of random variables. For such a reward, we
prove a new existence result for optimal stopping times under weaker
assumptions than in the classical case. This result is used to prove
the existence of optimal multiple stopping times for~$v(S)$ by
a~constructive method. Moreover, under strong regularity assumptions on~%
$\psi$, we show that the new reward $\phi$ can be aggregated by a
progressive process. This leads to new applications, particularly in
finance (applications to American options with multiple exercise
times).
\end{abstract}

%
\begin{keyword}[class=AMS]
\kwd[Primary ]{60G40}
\kwd[; secondary ]{60G07, 28B20, 62L15}.
\end{keyword}
\begin{keyword}
\kwd{Optimal stopping}
\kwd{optimal multiple stopping}
\kwd{aggregation}
\kwd{swing options}
\kwd{American options}.
\end{keyword}

\end{frontmatter}

\section*{Introduction}\label{in}

The present work on the optimal multiple stopping time problem,
following the optimal single stopping time problem, involves proving the
existence of the maximal reward, finding necessary or sufficient
conditions for the existence of optimal stopping times and providing a
method to compute these optimal stopping times.

The results are well known in the case of the optimal single stopping time
problem. Consider a \textit{reward} given by a right-continuous
left-limited (RCLL) positive adapted process $(\phi_t)_{{ 0\leq t\leq
T}}$ on
${\mathbb{F}}=(\Omega, \mathcal{F},(\mathcal{F}_t)_{{ 0\leq t\leq
T}},P)$, $\mathbb{F}$
satisfying the
usual conditions, and look for the maximal reward
\[
v(0)= \sup\{E[ \phi_\tau], \tau\in T_0\},
\]
where $T$ $\in$ $ ]0, \infty[ $ is the fixed time horizon and $T_0$
is the set of stopping times $\theta$ smaller than $T$. From now on,
the process $(\phi_t)_{{ 0\leq t\leq T}}$ will be denoted by $(\phi
_t)$. In order
to compute $v(0),$ we introduce for each $S\in T_0$ the \textit{value
function} $v(S)=\operatorname{ess}\sup\{ E[\phi_\tau |
\mathcal{F}_S], \tau\in T_S\}$,
where $T_S$ is the set of stopping times in $T_0$ greater than $S$. The
value function is given by a family of random variables $\{v(S), S\in
T_0\}$. By using the right continuity of the reward $(\phi_t)$, it
can be shown that there exists an adapted process $(v_t)$ which \textit{aggregates} the family of random variables $\{v(S), S\in T_0\}$
that is such that $v_S=v(S)$ a.s.~for each $S\in T_0$. This process is
the \textit{Snell envelope} of $(\phi_t)$, that is, the smallest
supermartingale process that dominates $\phi$. Moreover, when the
reward $(\phi_t)$ is continuous, the stopping time defined
trajectorially by
\[
\overline\theta(S)=\inf\{t\geq S, v_t=\phi_t\}
\]
is optimal.
For details, see El Karoui (\citeyear{EK}), Karatzas and Shreve (\citeyear{KS2}) or Peskir
and Shiryaev (\citeyear{Peskir}).

In the present work, we show that computing the value function for the
optimal multiple stopping time problem
\[
v(S)=\operatorname{ess}\sup\{E[\psi( \tau_1,\ldots,\tau_d) |
\mathcal{F}_S], \tau
_1,\ldots,\tau_d\in T_S\},
\]
reduces to computing the value function for an optimal single stopping
time problem
\[
u(S)= \operatorname{ess}\sup\{E[\phi(\theta) |\mathcal
{F}_S], \theta\in T_S \},
\]
where the \textit{new reward} $\phi$ is no longer an RCLL process, but a
family $\{\phi(\theta), \theta\in T_0\}$ of positive random
variables which satisfies some compatibility properties. For this new
optimal single stopping time problem with a reward $\{\phi(\theta),
\theta\in T_0\}$, we show that the minimal optimal stopping time for
the value function $u(S)$ is no longer given by a hitting time of
processes, but by
the essential infimum
\[
\theta^*(S) :=\operatorname{ess}\inf\{
\theta\in T_S
,  u(\theta) =\phi(\theta) \mbox{ a.s.} \}.
\]
This method also has the advantage that it no longer requires any
aggregation results that need stronger hypotheses and whose proofs are
rather technical.

By using the reduction property $v(S) = u(S)$ a.s., we give a method to
construct by induction optimal stopping times $(
\tau^*_1,\ldots,\tau^*_d)$ for $v(S)$, which are also defined as
essential infima, in terms of \textit{nested} optimal single stopping time
problems.

Some examples of optimal multiple stopping time problems have been
studied in different mathematical fields. In finance, this type of
problem appears in, for instance, the study of \textit{swing options}
[e.g., Carmona and Touzi (\citeyear{CT}), Carmona and Dayanik (\citeyear{CD})] in the
case of ordered stopping times. In the nonordered case, some optimal
multiple stopping time problems appear as useful mathematical tools to
establish some large deviations estimations [see Kobylanski and Rouy
(\citeyear{KR})]. Further applications can be imagined in, for example, finance
and insurance [see Kobylanski, Quenez and Rouy-Mironescu (\citeyear{KQR})]. In a work in preparation
[see Kobylanski and Quenez (\citeyear{KQ})], the Markovian case will be studied
in detail and some applications will be presented.

The paper is organized as follows. In Section \ref{1} we revisit the optimal
single stopping time problem for admissible families. We prove the
existence of optimal stopping times when the family $\phi$ is right-
and left-continuous in expectation along stopping times. We also
characterize the minimal optimal stopping times. In Section \ref{2} we solve
the optimal double stopping time problem. Under quite weak assumptions,
we show the existence of a pair of optimal stopping times and give a
construction of those optimal stopping times. In Section \ref{3} we
generalize the results obtained in Section \ref{2} to the optimal
$d$-stopping-times problem. Also, we study the simpler case of a
symmetric reward.
In this case, the problem clearly reduces to ordered stopping times,
and our general characterization of the optimal multiple stopping time
problem in terms of nested optimal single stopping time problems
straightforwardly reduces to a sequence of optimal single stopping time
problems defined by backward induction.
We apply these results to \textit{swing options} and, in this particular
case, our results correspond to those of Carmona and Dayanik (\citeyear{CD}).
In the last section, we prove some aggregation results and characterize
the optimal stopping times in terms of hitting times of processes.

Let ${\mathbb F}=(\Omega, \mathcal{F}, (\mathcal{F}_t)_{{ 0\leq
t\leq T}},P)$ be a probability space,
where $T\in\,]0, \infty[$ is the fixed time horizon and $(\mathcal
{F}_t)_{{ 0\leq t\leq T}}$
is a filtration satisfying the usual conditions of right continuity and
augmentation by the null sets of $\mathcal{F}= \mathcal{F}_T$. We suppose
that $\mathcal{F}_0$ contains only sets of probability $0$ or $1$. We
denote by~$T_{0}$ the collection of stopping times of ${\mathbb{F} }$
with values in $[0 , T]$. More generally, for any stopping time $S$,
we denote by $T_{S}$ the class of stopping times $\theta\in T_0$ with
$S \leq\theta$ a.s.

We use the following notation: for real-valued random variables $X$ and
$X_n$, $n \in\mathbb{N}$, the notation ``$X_n\uparrow X$'' means ``the
sequence $(X_n)$ is nondecreasing and converges to $X$ a.s.''

\section{The optimal single stopping time problem revisited}\label{1}

We first recall some classical results on the optimal single stopping time
problem.

\subsection{Classical results}\label{1.1}

The following classical results, namely the supermartingale property of
the value function, the optimality criterium and the right continuity
in expectation of the value function are well known [see El Karoui
(\citeyear{EK}) or Karatzas and Shreve (\citeyear{KS2}) or Peskir and Shiryaev (\citeyear{Peskir})].
They are very important tools in optimal stopping theory and will
often be used in this paper in the (unusual) case of a reward given by
an admissible family of random variables defined as follows.

\begin{Definition}
A family of random variables $\{\phi(\theta), \theta\in T_0\}$ is
said to be
\textit{admissible} if it satisfies the following conditions:\vspace*{-2pt}
\begin{enumerate}
\item
for all
$\theta\in T_0$, $\phi(\theta)$ is an $\mathcal{F}_\theta
$-measurable $\overline {\mathbb{R}}^+$-valued
random variable;
\item
for all
$\theta,\theta'\in T_0$, $\phi(\theta)=\phi(\theta')$ a.s.~on
$\{\theta=\theta'\}$.\vspace*{-2pt}
\end{enumerate}
\vspace*{-2pt}
\end{Definition}

\begin{Remark}
Let $(\phi_t )$ be a positive progressive process. The family defined by
$\phi(\theta)=\phi_\theta$ is admissible.

Note also that the definition of admissible families corresponds to the
notion of $T_0$-systems introduced by El Karoui (\citeyear{EK}).\vspace*{-2pt}
\end{Remark}

For the convenience of the reader, we recall the definition of the
essential supremum and its main properties in Appendix
\hyperref[app1]{A}.

Suppose the reward is given by an admissible family
$\{\phi(\theta),\theta\in T_0\}$.
The \textit{value function at time $S$}, where $S\in T_0$, is given
by\vspace*{-1pt}
%
\begin{equation}\label{vs}
v(S)=\operatorname{ess}\sup_{\theta\in T_S} E[\phi(\theta)  |
\mathcal{F}_S].\vspace*{-2pt}
\end{equation}

\begin{Proposition}[(Admissibility of the value function)]\label{P1.Adm}
The value function that is the family of random variables $\{v(S),
S\in T_0\}$ defined by (\ref{vs}) is an admissible family.\vspace*{-2pt}
\end{Proposition}

\begin{pf}
Property 1 of admissibility for $\{v(S), S\in T_0\}$ follows from the
existence of the essential supremum (see Theorem \ref{TA} in Appendix
\hyperref[app1]{A}).

Take $S,S'\in T_0$ and let $A=\{S=S'\}$. For each $\theta\in T_S$, put
$\theta_A=\theta{\mathbf{1}}_A+T {\mathbf{1}}_{A^c}$. As $A \in\mathcal
{F}_S\cap
\mathcal{F}_{S'}$, we have a.s.~on $A$, $E[\phi(\theta)|\mathcal{F}_S]=
E[\phi(\theta_A)| \mathcal{F}_S]=$ $E[\phi(\theta_A)|
\mathcal{F}_{S'}]\leq
$ $
v(S'),$ hence taking the essential supremum over $\theta\in T_S$, we
have $v(S)\leq v(S')$ a.s., and by symmetry of $S$ and $S'$, we have
shown property 2 of admissibility.\vspace*{-2pt}
\end{pf}

\begin{Proposition}\label{P1.2a}
There exists a sequence of stopping times $(\theta^n)_{n \in\mathbb
{N}}$ with
$\theta^n $ in $ T_{S}$ such that\vspace*{-2pt}
\[
E[\phi(\theta^{n})  |\mathcal{F}_{S}]\ \uparrow\ v(S)\qquad\mbox{a.s.}\vspace*{-2pt}
\]
\end{Proposition}

\begin{pf}
For each $S$ $\in$ $T_0$, one can show that the
set $\{ E[\phi(\theta) |\mathcal{F}_S],  \theta\in T_S \}
$ is closed
under pairwise maximization. Indeed, let $\theta,\theta'\in T_0$ and
$A=\{ E[\phi(\theta')| \mathcal{F}_S]\leq
E[\phi(\theta)| \mathcal{F}_S]\}$. One has $A\in\mathcal
{F}_S$. Let $\tau
=\theta{\mathbf{1}}_A+\theta' {\mathbf{1}}_{A^c}$, a stopping time. It is easy
to check that $ E[\phi(\tau)| \mathcal{F}_S]=$ $ E[\phi(\theta
)| \mathcal{F}
_S]\vee$ $ E[\phi(\theta')| \mathcal{F}_S]$.
The result follows by a classical result (see Theorem \ref{TA} in
Appendix \hyperref[app1]{A}).\vspace*{-2pt}
\end{pf}

Recall that for each fixed $S\in T_0$, an admissible family $\{
h(\theta) ,
 \theta\in T_{S} \}$ is said to be a \textit{supermartingale system}
(resp., a \textit{martingale system}) if, for any~$\theta$, $\theta'$ $\in$ $T_0$
such that $\theta\geq\theta'$ a.s.,\vspace*{-1pt}
\[
E[h(\theta)  |\mathcal{F}_{ \theta' }] \leq h(\theta')
\qquad\mbox{a.s.}\qquad\bigl(\mathrm{resp.,}~E[h(\theta)  |\mathcal{F}_{ \theta'
}] = h(\theta') \mbox{ a.s.}\bigr).
\]\vadjust{\eject}

\begin{Proposition}\label{P1.SuperM}
\begin{itemize}
\item The value function $\{v(S) ,  S \in T_0 \}$ is a
supermartingale system.
\item Furthermore, it is characterized as the
\textup{Snell envelope system} associated with $\{ \phi(S) , S\in T_0
\}$,
that is, the smallest supermartingale system which is greater (a.s.)
than $\{ \phi(S) , S\in T_0\}$.
\end{itemize}
\end{Proposition}

\begin{pf}
Let us prove the first part. Fix $S \geq S'$
a.s. By Proposition \ref{P1.2a}, there exists an optimizing sequence
$(\theta^n)$ for $v(S)$. By the monotone convergence theorem,
${ E[v(S)  | \mathcal{F}_{ S' }] =
\lim_{n\to\infty} E[\phi(\theta^{n})  |\mathcal{F}_{ S' }]}$ a.s.
Now, for each $n$, since $\theta^n\geq S'$ a.s., we have
$E[\phi(\theta^n)|\mathcal{F}_{S'}]\leq v(S')$ a.s. Hence,
$E[v(S)  |
\mathcal{F}_{ S' }] \leq v(S')$ a.s., which gives
the supermartingale property of the value function.

Let us prove the second part. Let $\{v'(S),S\in T_0\}$ be a supermartingale
system such that for each $\theta\in T_0$, $v'(\theta)$ $\geq$
$\phi(\theta)$ a.s. Fix $S\in T_0$. By the properties of $v'$, for
all $\theta\in T_S$,
$v'(S)\geq E[v'(\theta)| \mathcal{F}_S]\geq$ $E[\phi(\theta)
| \mathcal{F}_S]$
a.s. Taking the supremum over $\theta\in T_S$, we have $v'(S)\geq
v(S)$ a.s.
\end{pf}

Now, recall the following Bellman optimality criterium [see, e.g., El
Karoui (\citeyear{EK})].

\begin{Proposition}[(Optimality criterium)]\label{criterium}
Fix $S$ $\in$ $T_0$ and let $\theta^* \in T_S$ be such that
$E[\phi(\theta^*)]<\infty$. The three following assertions are
equivalent:
\begin{enumerate}
\item$\theta^*$ is $S$-optimal for $v(S)$, that is,
%
\begin{equation}\label{So}
v(S) = E[\phi( \theta^*)  |\mathcal{F}_{S}] \qquad\mbox{a.s.};
\end{equation}
\item$v ( \theta^*) = \phi( \theta^*) \mbox{ a.s. and }
 E[v(S)]= E[ v( \theta^*)];$
\item$E[v(S)]= E[ \phi( \theta^*)].$
\end{enumerate}
\end{Proposition}

\begin{Remark}
Note that since the value function is a supermartingale
system, equality $E[v(S)]= E[ v( \theta^*)]$ is equivalent to the fact
that the family $\{ v(\theta), \theta\in T_{S, \theta^*} \}$ is
a martingale system.
\end{Remark}

\begin{pf*}{Proof of Proposition \ref{criterium}} Let us show that assertion 1 implies
assertion~2. Suppose assertion 1 is satisfied. Since the value function
$v$ is a~supermartingale system greater that $\phi$, we clearly have
\[
v(S) \geq E[v(\theta^*)  |\mathcal{F}_{ S }] \geq E[\phi(\theta^*)
 |\mathcal{F}_{ S }]\qquad \mbox{a.s. }
\]
Since equality (\ref{So}) holds, this implies that the previous
inequalities are actually equalities.

In particular, $E[v(\theta^*)  |\mathcal{F}_{ S }] = E[\phi
(\theta^*)  |\mathcal{F}_{ S }] $ a.s., but as inequality
$v(\theta^*) \geq\phi(\theta^*)$ holds a.s., and as $E[\phi(\theta
^*)]<\infty$, we have $v(\theta^*) = \phi(\theta^*)$ a.s.

Moreover, $v(S) = E[v(\theta^*)  |\mathcal{F}_{ S }] $ a.s., which
gives $E[v(S)]= E[ v( \theta^*)]$. Hence, assertion 2 is satisfied.

Clearly, assertion 2 implies assertion 3. It remains to show that 3
implies~1.

Suppose that 3 is satisfied. Since $v(S) \geq E[\phi(\theta^*)
|\mathcal{F}_{ S }] $ a.s., this gives $v(S) = E[\phi(\theta^*)
|\mathcal{F}_{ S }] $ a.s. Hence, 1 is satisfied.
\end{pf*}

\begin{Remark} \label{eo}
It is clear that
%
\begin{equation}\label{nstp}
E[v(S)]= \sup_{\theta\in T_S} E[ \phi( \theta)].
\end{equation}
By assertion 3 of Proposition \ref{criterium}, a stopping time
$\theta^* \!\!\in T_S$ such that $E[\phi(\theta^*)]\!<\infty$ is $S$-optimal
for $v(S)$ if and only if it is optimal for the optimal stopping time
problem (\ref{nstp}), that is,
\[
\sup_{\theta\in T_S} E[ \phi( \theta)]= E[\phi(\theta^*)].
\]
%

\end{Remark}

We now give a regularity result on $v$ [see Lemma 2.13 in El Karoui
(\citeyear{EK})]. Let us first introduce the following definition.

\begin{Definition}
An admissible family $\{\phi(\theta),\theta\in T_0\}$
is said to
be \textit{right- (resp., left-) continuous along stopping times in
expectation [RCE (resp., LCE)]} if for any $\theta\in T_0$ and any
sequence $(\theta_n)_{n\in\mathbb{N}}$ of stopping times such that
$\theta_n\downarrow\theta$ a.s.~(resp.,~$\theta_n\uparrow\theta$ a.s.),
one has $ {E[\phi(\theta)]=\lim_{n\to\infty}
E[\phi(\theta_n)].} $
\end{Definition}

\begin{Remark}
If $(\phi_t)$ is a continuous adapted process such that\break
$E[{\sup_{t\in[0,T]}\phi_t] < \infty}$, then the family
defined by $\phi(\theta)=\phi_\theta$ is clearly RCE and LCE. Also, if
$(\phi_t )$ is an RCLL adapted process such that its jumps are totally
inaccessible, then the family
defined by $\phi(\theta)=\phi_\theta$ is clearly RCE and even LCE.
\end{Remark}

\begin{Proposition}\label{L1}
Let $\{\phi(\theta),\theta\in T_0\}$ be an admissible family
which is RCE.
The family $\{ v(S) ,  S \in T_0 \}$ is then RCE.
\end{Proposition}

\begin{pf}
Since $\{v(S),S\in T_0\} $ is a supermartingale
system, the function $S\mapsto E[ v(S)]$ is a nonincreasing function of
stopping times.
Suppose it is not RCE at $S\in T_0$. If $E[v(S) ]< \infty$, then there
exists a constant $\alpha>0$ and a~%
sequence of stopping times $(S_n)_{n\in\mathbb{N}}$ such that
$S_n\downarrow
S$ a.s.~and such that
%
\begin{equation}\label{vSna}\lim_{n\to\infty}\uparrow
E[v(S_n)]+\alpha\leq E[v(S)].
\end{equation}
Now, recall that ${E[v(S)]= \sup_{\theta\in T_S}E[\phi
(\theta
)]}$ [see (\ref{nstp})]. Hence,
there exists $\theta'\in T_S$ such that
%
\[
\sup_{n \in\mathbb{N}} \sup_{\theta\in T_{S_n}}E[\phi(\theta
)]+\frac
{\alpha}{2} \leq E[\phi(\theta')].
\]
Hence, for all $n\in\mathbb{N}$,
$ E[\phi(\theta'\vee S_n)]+\frac{\alpha}{2}\leq E[\phi
(\theta')]$.
As $\theta'\vee S_n \downarrow\theta'$ a.s., we~ob\-tain, by taking
the limit when $n\to\infty$ and using the RCE property of $\phi$,~that
\[
E[\phi(\theta')]+\frac{\alpha}{2} \leq E[\phi(\theta')],
\]
which gives the expected contradiction in the case $E[v(S) ]< \infty
$.

Otherwise, instead of (\ref{vSna}), we have
${\lim_{n\to\infty}\uparrow E[v(S_n)] \leq C}$ for
some constant $C>0$, and similar arguments as in the finite case lead
to a contradiction as well.
\end{pf}

\subsection{New results}\label{1.2}

We will now give a new result which generalizes the classical existence
result of an optimal stopping time stated in the case of a~reward
process to the case of a reward family of random variables.

\begin{Theorem}[(Existence of optimal stopping times)]\label{T.1}
Let $\{ \phi(\theta), \theta\in T_0\}$ be an admissible family that
satisfies the integrability condition
\[
v(0) = \sup_{\theta\in T_0} E[  \phi(\theta)]<\infty
\]
 and which is RCE and LCE along stopping times.
Then, for each $S$ $\in$ $T_0$, there exists an optimal stopping time
for $v(S)$. Moreover, the random variable defined by
%
\begin{equation}\label{E.te}\theta^*(S) :=\operatorname{ess}\inf\{
\theta\in
T_S ,  v(\theta) =\phi(\theta) \mbox{ a.s.} \}
\end{equation}
is the minimal optimal stopping time for $v(S)$.
\end{Theorem}

Let us emphasize that in this theorem, the optimal stopping time
$\theta^*(S)$ is not defined trajectorially, but as an essential
infimum of random variables. In the classical case, that is, when the
reward is given by an adapted RCLL process, recall that the minimal
optimal stopping time is given by the random variable $\overline\theta
(S)$ defined trajectorially by
\[
\overline\theta(S)=\inf\{t\geq S, v_t=\phi_t\}.
\]

The definition of $\theta^*(S)$ as an essential infimum allows the
assumption on the regularity of the reward to be relaxed. More
precisely, whereas in the previous works (mentioned in the
\hyperref[in]{Introduction}), the reward was given by an RCLL and LCE process, in our
setting, the reward is given by an RCE and LCE family of random
variables. The idea of the proof is classical: we use an approximation
method introduced by Maingueneau (\citeyear{Ma}), but our setting allows us to
simplify and shorten the proof.

\begin{pf*}{Proof of Theorem \ref{T.1}} The proof will be divided into two parts.

\textit{Part I}: In this part, we will prove the existence of an optimal
stopping time.

Fix $S$ $\in$ $T_0$. We begin by constructing a family of stopping
times [see~Main\-gueneau (\citeyear{Ma}) or El Karoui (\citeyear{EK})]. For $\lambda$ $\in$
$]0,1[$, define the $\mathcal{F}_S$-measurable random variable $\theta
^{\lambda}
(S)$ by
%
\begin{equation} \label{tlS}
\theta^{\lambda} (S) := \operatorname{ess}\inf\{\theta\in T_S
, \lambda
v(\theta) \leq\phi(\theta) \mbox{ a.s.} \}.
\end{equation}

The following lemma holds.

\begin{Lemma}\label{eli} The stopping time
$\theta^{\lambda} (S)$ is a $(1- \lambda)$-optimal stopping time for
%
\begin{equation} \label{oie}
 E[v(S)]= \sup_{\theta\in T_S} E[ \phi( \theta)],
\end{equation}
that is,
%
\begin{equation}\label{la}
\lambda E[v(S)]  \leq E [ \phi(\theta^{\lambda}(S) ) ].
\end{equation}
\end{Lemma}

Suppose now that we have proven Lemma \ref{eli}.

Since $\lambda\mapsto\theta^\lambda(S)$ is nondecreasing, for $S\in
T_0$, the stopping
time
%
\begin{equation}\label{chapeau}
\hat\theta(S):= \lim_{\lambda\uparrow1} \uparrow\theta^\lambda
(S)
\end{equation}
is well defined. Let us show that $\hat{\theta}(S)$
is optimal for $v(S)$.


By letting $\lambda\uparrow1$ in inequality (\ref{la}), and since
$\phi$ is LCE, we easily derive that $ E[v(S )]$ $ =$
$E[\phi(\hat{\theta}(S) ) ]$. Consequently, by the optimality criterium
3 of Proposition~\ref{criterium}, $\hat{\theta}(S)$ is $S$-optimal for
$v(S)$. This completes part I.

\textit{Part II}: Let us now prove that $\theta^*(S)=\hat
\theta(S)$ a.s., where $\theta^*(S)$ is defined by (\ref{E.te}), and
that it is the minimal optimal stopping time for $v(S)$.

For each $S$ $\in$ $T_0$, the set ${\mathbb T}_S= \{\theta\in T_S
,  v(\theta) =\phi(\theta) \mbox{ a.s.} \}$ is not empty (sin\-ce~%
$T$ belongs to ${\mathbb T}_S$) and is closed under pairwise
minimization. Hence, there exists a sequence $(\theta_n) _{n \in
\mathbb{N}}$
of stopping times in ${\mathbb T}_S$ such that $ \theta_n\downarrow
\theta^* (S)\mbox{ a.s.}$ Consequently, $\theta^* (S)$ is a
stopping time.

Let $\theta$ be an optimal stopping time for $v(S)$. By the optimality
criterium (Proposition \ref{criterium}), and since, by assumption,
$E[\phi(\theta)] < \infty$, we have \mbox{$v(\theta)=\phi(\theta)$} a.s.~and
hence
\[
\theta^*(S)\leq\operatorname{ess}\inf\{\theta\in T_0,  \theta
\mbox{ optimal for }
v(S)\}\qquad \mbox{a.s. }
\]

Now, for each $\lambda<1,$ the stopping time $\theta^\lambda(S)$
defined by
(\ref{tlS}) clearly satisfies $\theta^\lambda(S)\leq\theta^*(S)$
a.s. Passing to
the limit when $\lambda\uparrow1,$ we obtain $\hat\theta(S) \leq
\theta^*(S)$. As $\hat\theta(S)$ is optimal for $v(S)$, this implies
that $\hat\theta(S) \geq\operatorname{ess}\inf\{\theta\in T_0,
\theta$ optimal
for $
v(S)\} $ a.s. Hence,
\[
\theta^ *(S)=\hat\theta(S)= \operatorname{ess}\inf\{\theta\in
T_0,  \theta\mbox{
optimal for } v(S)\}\qquad \mbox{a.s.},
\]
which gives the desired result. The proof of Theorem \ref{T.1} is thus
complete.
\end{pf*}

 It now remains to prove Lemma \ref{eli}.

\begin{pf*}{Proof of Lemma \ref{eli}}
 We have to prove inequality
(\ref{la}). This will be done by means of the following steps.

\textit{Step 1}: Fix $\lambda\in\,]0,1[$. It is easy to
check that the set ${\mathbb{T}}^{\lambda}_{S}=\{\theta\in T_S ,
\lambda
v(\theta) \leq\phi(\theta) \mbox{ a.s.} \}$ is nonempty
(since $T$ $\in$ ${\mathbb{T}}^{\lambda}_{S}$) and closed by
pairwise minimization. By Theorem \ref{TA} in the \hyperref[app1]{Appendix}, there
exists a sequence $(\theta^n)$ in ${\mathbb{T}}_{S}$ such that
$\theta^n\downarrow\theta^\lambda(S)$ a.s. Therefore, $\theta
^\lambda(S)$ is a
stopping time and $\theta^\lambda(S)\geq S$ a.s. Moreover, we have
$\lambda
v(\theta^n)\leq\phi(\theta^n)$ a.s.~for all $n$. Taking
expectation and using the RCE properties of $v$ and $\phi$, we obtain
%
\begin{equation}\label{E.lvphi}
\lambda E[ v(\theta^{\lambda} (S)] ) \leq E[\phi(\theta^{\lambda}
(S) )
] .
\end{equation}

\textit{Step 2}: Let us show that for each
$\lambda\in\,]0,1[$ and each $S\in T_0$,
%
\begin{equation}\label{E.u}
v(S) = E [ v(\theta^{\lambda} (S)) |  \mathcal{F}_S]\qquad\mbox{a.s.}
\end{equation}
For each $S$ $\in$ $T_0$, let us define the random variable $J(S) = E
[ v(\theta^{\lambda} (S))  | \mathcal{F}_S].$ Step 2 amounts to
showing that $J(S)=v(S)$ a.s.

Since $\{ v(S), S \in T_0 \}$ is a supermartingale system
and since $\theta^{\lambda} (S) \geq S$ a.s., we have that
\[
J(S) = E [ v(\theta^{\lambda} (S))  | \mathcal{F}_S] \leq v(S)
\qquad\mbox{a.s.}
\]

It remains to show the reverse inequality.

\textit{Step 2a}: Let us show that the family $\{ J(S),
S\in T_0\}$ is a supermartingale system.

Let $S, S' \in T_0 $ be such that $S' \geq S$ a.s. As
$\theta^{\lambda} (S')  \geq \theta^{\lambda} (S) \geq S$
a.s., we have
\[
E[ J(S' )  | \mathcal{F}_S] = E [ v(\theta^{\lambda} (S'))
| \mathcal{F}_S]
= E\bigl[ E \bigl[ v(\theta^{\lambda} (S'))  | \mathcal{F}_{\theta
^{\lambda} (S)} \bigr]  |\mathcal{F}_S \bigr]\qquad\mbox{a.s.}
\]
Now, since $\{ v(S), S\in T_0\}$ is a supermartingale system, $E [
v(\theta^{\lambda} (S'))  | \mathcal{F}_{\theta^{\lambda} (S)}
] $
$ \leq v({\theta^{\lambda} (S)}) \mbox{ a.s.}$ Consequently,
\[
E[ J(S' )  | \mathcal{F}_S] \leq E [ v(\theta^{\lambda} (S))
| \mathcal{F}_S]= J(S)\qquad \mbox{a.s.}
\]

\textit{Step 2b}: Let us show that for each $ S \in T_0$
and each $\lambda\in\,]0,1[$,
\[
\lambda v(S) + (1 - \lambda) J(S)  \geq\phi(S)\qquad \mbox{a.s.}
\]
Fix $ S \in T_0$ and $\lambda\in\,]0,1[$.

On $\{ \lambda v(S) \leq\phi(S) \}$, we have $ \theta^{\lambda}
(S) = S $ a.s. Hence, on $\{ \lambda v(S) \leq\phi(S) \}$, $J(S)
= E [ v(\theta^{\lambda} (S))  | \mathcal{F}_S] = E [ v(S)
 | \mathcal{F}_S]  =  v(S) $ and therefore
\[
\lambda v(S) + (1 -
\lambda) J(S) = v(S) \geq\phi(S)\qquad\mbox{a.s.}
\]
Furthermore,
on $\{ \lambda v(S) > \phi(S) \}$, as $J(S)$ is nonnegative, we have
\[
\lambda v(S) + (1 - \lambda) J(S)  \geq \lambda v(S) \geq
\phi(S)\qquad \mbox{a.s.},
\]
and the proof of Step 2b is complete.

Now, the family $ \{\lambda v(S) + (1 - \lambda) J(S) ,S\in T_0\}$
is a
supermartingale system by convex combination of two supermartingale
systems. Hence, as the value function $\{ v(S),S\in T_0\}$ is
characterized as the smallest supermartingale system which dominates
$\{ \phi(S), S\in T_0\}$, we derive that for each $S \in T_0$,
\[
\lambda v(S) + (1 - \lambda) J(S) \geq v(S)\qquad \mbox{a.s.}
\]
Now, by the integrability assumption made on $\phi$, we have $v(S) <
\infty$ a.s. Hence, we have $J(S) \geq v(S)$ a.s. Consequently, for
each $S \in T_0$, $J(S) = v(S)$ a.s., which completes Step 2.

Finally, Step 1 [inequality (\ref{E.lvphi})] and Step 2 [equality
(\ref{E.u})] give
\[
\lambda E[v(S)] =  \lambda E [ v(\theta^{\lambda}(S) )]  \leq
 E [ \phi(\theta^{\lambda}(S) ) ].
\]
In other words, $\theta^{\lambda} (S)$ is a $(1- \lambda)$-optimal
stopping time for (\ref{oie}), which completes the proof of Lemma
\ref{eli}.
\end{pf*}

\begin{Remark}
Recall that in the previous works [see, e.g., Karatzas and Shreve
(\citeyear{KS2}), Proposition D.10 and Theorem D.12], the proof of the existence
of optimal stopping times requires the value function to be aggregated
and thus the use of some fine aggregation results such as Proposition
\ref{P.SMA}. In our work, since we only work with families of random
variables, we do not need any aggregation techniques, which
simplifies and shortens the proof.
\end{Remark}

Under some regularity assumptions on the reward, we can show that the
value function family is left-continuous along stopping times in expectation. More
precisely, we have the following.


\begin{Proposition}\label{P.vLCE}
$\!\!\!\!\!$Suppose that the admissible family $\{ \phi(\theta),
\theta\!\in\!T_0\}$ is~LCE
and RCE, and satisfies the integrability condition \mbox{$v(0) = \sup_{\theta\in
T_0} E[\phi(\theta)]<\infty$}.

The value function $\{v(S), S\in T_0\}$ defined by (\ref{vs}) is then
LCE.
\end{Proposition}

\begin{pf}
Let $S\in T_0$ and let $(S_n)$ be a sequence of
stopping times such that $S_n\uparrow S$ a.s. Note that by the
supermartingale property of $v,$ we have
%
\begin{equation}\label{E.vSn}
E[v(S_n)]\geq E[v(S)].
\end{equation}
Now, by Theorem \ref{T.1}, the stopping time $\theta^*(S_n)$ defined by
(\ref{E.te}) is optimal for $v(S_n)$. Moreover, it is clear that
$(\theta^*(S_n))_n$ is a nondecreasing sequence of stopping times
dominated by $\theta^*(S)$.

Let us define ${\overline\theta=
\lim_{n\to\infty} \uparrow\theta^*(S_n)}$. Note that $\overline
\theta$ is a stopping time. Also, as for each $n$, $\theta^*(S_n)\geq
S_n$ a.s., it follows that $\overline\theta\geq S$ a.s. Therefore,
since $\phi$ is LCE,
\[
E[ v(S)]\geq E[\phi(\overline\theta)]= \lim_{n\to\infty} E[\phi
(\theta^*(S_n))]=\lim_{n\to\infty} E[v(S_n)].
\]
This, together with (\ref{E.vSn}), gives ${ E[
v(S)]=\lim_{n\to\infty} E[v(S_n)]}$.
\end{pf}

\begin{Remark}\label{unun}
$\!\!\!$In this proof, we have also proven that $\overline\theta$ is
optimal~for~$v(S)$. Hence, by the optimality criterium, $v(\overline\theta)=
\phi(\overline\theta)$ a.s., which implies that $\overline\theta
\geq\theta^*(S)$ a.s. Moreover, since for each $n$, $\theta^*(S_n)
\leq\theta^*(S)$ a.s., by letting~$n$ tend to $\infty$, we clearly
have that $\overline\theta\leq\theta^*(S)$ a.s. Hence,
${\overline\theta= \lim_{n\to\infty} \uparrow
\theta^*(S_n)}$ $= $ $ \theta^*(S)$ a.s. Thus, we have also shown that
the map $S \mapsto\theta^* (S)$ is left-con\-tinuous along stopping
times.
\end{Remark}

\section{The optimal double stopping time problem}\label{2}
\subsection{Definition and first properties of the value function}\label{2.1}

We now consider the optimal double stopping time problem. We introduce
the following definitions.\vadjust{\goodbreak}

\begin{Definition} The family
$\{\psi(\theta,S),\theta,S\in T_0\}$ is \textit{biadmissible} if it
satisfies:
\begin{enumerate}
\item
for all $ \theta,S\in T_0,$ $\psi(\theta,S)$ is an $\mathcal
{F}_{\theta\vee
S}$-measurable $\overline {\mathbb{R}}^+$-valued r.v.;
\item
for all $\theta,\theta',S,S'\in
T_0,  \psi(\theta,S)=\psi(\theta',S')$ a.s.~on $
\{\theta=\theta'\}\cap\{S=S'\} $.
\end{enumerate}
\end{Definition}

\begin{Remark}\label{biprocess}
Let $\Psi$ be a biprocess, that is, a function
\[
\Psi\dvtx [0,T] ^2\times\Omega\rightarrow\mathbb{R}^+;
(t,s,\omega)
\mapsto\Psi_{t,s} (\omega)
\]
such that for almost all $\omega$, the map $(t,s) \mapsto\Psi_{t,s}
(\omega)$ is \textit{right-continuous} (i.e., $\Psi_{t,s}$ $=$
${ \lim_{ (t' , s') \to(t^+, s^+)} \Psi_{t',s'}}$), and
for each $(t,s) \in[0,T] ^2$, $\Psi_{t,s}$ is $\mathcal{F}_{t \vee
s}$-measurable.
In this case, the family $\{\psi(\theta,S),\theta,S\in T_0\}$
defined by
\[
\psi(\theta,S)(\omega) := \Psi_{\theta(\omega), S(\omega)}
(\omega)
\]
is clearly biadmissible.
\end{Remark}

For a biadmissible family $\{\psi(\theta,S),\theta,S\in T_0\},$
let us consider
the value function associated with the reward family
$\{\psi(\theta,S),\theta,S\in T_0\}$:
%
\begin{equation}\label{vS}
v(S)=\operatorname{ess}\sup_{\tau_1,\tau_2\in T_S} E[\psi(\tau
_1,\tau_2) |\mathcal{F}_S].
\end{equation}

As in the case of the single stopping time problem, we have the following
properties.

\begin{Proposition}\label{P2.1}
Let $\{\psi(\theta,S),\theta,S\in T_0\}$ be a biadmissible
family of random
variables. The following properties then hold:
\begin{enumerate}[(1)]
\item[(1)] the family $\{v(S), S\in T_0\}$ is an admissible family of random
variables;
\item[(2)] for each $S\in T_0$, there exists a sequence of pairs of
stopping times $((\tau^n _1,\break \tau^n _2))_{n\in\mathbb
{N}}$ in
$ T_S
\times T_S$ such that $\{  E[\psi(\tau^n_1,\tau^n_2) |\mathcal
{F}_S]
\}_{n \in\mathbb{N}}$ is nondecreasing and a.s.
\[
E[\psi(\tau^n_1,\tau^n_2) |\mathcal{F}_S] \uparrow v(S);
\]
\item[(3)]the family of random variables $\{ v(S) , S\in T_0\}$ is a
supermartingale system, that is, it satisfies the dynamic programming
principle.
\end{enumerate}
\end{Proposition}

\begin{pf} (1) As in the case of single stopping time,
property 1 of admissibility for $\{v(S), S\in T_0\}$ follows from the
existence of the essential supremum.

Take $S,S'\in T_0$ and put $A=\{S=S'\}$, and for each
$\tau_1,\tau_2\in T_S,$ put $\tau_1^A=\tau_1 {\mathbf{1}}_A+T
{\mathbf{1}}_{A^c}$ and $\tau_2^A=\tau_2 {\mathbf{1}}_A+T {\mathbf{1}}_{A^c}$. As
$A \in\mathcal{F}_S\cap\mathcal{F}_{S'}$, one has, a.s.~on~$A$,
$E[\psi(\tau_1,\tau_2)|\mathcal{F}_S]= E[\psi(\tau_1^A,\tau
_2^A)|
\mathcal{F}_S]=$ $E[\psi(\tau_1^A,\tau_2^A)| \mathcal
{F}_{S'}]\leq$ $ v(S').$ Hence,
taking the essential supremum over $\tau_1,\tau_2\in T_S$, we have
$v(S)\leq v(S')$ a.s., and, by symmetry, we have shown property 2 of
admissibility. Hence, the family $\{v(S), S\in T_0\}$ is an admissible
family of random variables.

The proofs of (2) and (3) can be easily adapted from the proofs of
Proposition~\ref{P1.2a} and Proposition~\ref{P1.SuperM}.
\end{pf}




Following the case of single stopping time, we now give some regularity
results on the value function.

\begin{Definition}\label{RCE2}
A biadmissible family $\{\psi(\theta,S),\theta,S\in T_0\}$ is
said to be
\textit{right-continuous along stopping times in expectation (RCE)} if,
for any~$\theta$, $S\in T_0$ and any sequences $(\theta_n)_{n\in\mathbb
{N}}$ $\in$ $T_0$
and $(S_n)_{n\in\mathbb{N}}$ $\in$ $T_0$ such that $\theta
_n\downarrow
\theta$
and $S_n \downarrow S$ a.s., one has $
E[\psi(\theta,S)]=\lim_{n\to\infty} E[\psi(\theta_n,S_n)]. $
\end{Definition}

\begin{Proposition}\label{L2.vRCE}
Suppose that the biadmissible family $\{\psi(\theta,S), \theta
,S\in T_0\}$
is RCE. The family $\{ v(S) , S\in T_0\}$ defined by (\ref{vS}) is
then RCE.
\end{Proposition}

\begin{pf}
The proof follows the proof of Proposition
\ref{L1}.
\end{pf}

\subsection{Reduction to an optimal single stopping time problem}\label{2.2}

In this section, we will show that the optimal double stopping time
problem (\ref{vS}) can be reduced to an optimal single stopping time
problem associated with a~\textit{new reward family}.

More precisely, for each stopping
time $\theta\in T_S$ let us introduce the two $\mathcal{F}_{\theta
}$-measurable
random variables
%
\begin{equation}\label{u12}
\qquad u_1({\theta})=\operatorname{ess}\sup_{\tau_1\in T_{\theta}} E[\psi(\tau_1,
\theta) |\mathcal{F}_{\theta}],\qquad u_2({\theta}) =\operatorname{ess}\sup_{\tau_2\in T_{\theta}} E[\psi(\theta,\tau_2)
|\mathcal{F}_{\theta}].\hspace*{-5pt}
\end{equation}
Note that since $\{\psi(\theta,S),\theta,S\in T_0\} $ is
biadmissible, for
each fixed $\theta\in T_0$, the families $\{\psi(\tau_1,\theta),
\tau_1 \in T_0\}$ and $\{\psi(\theta, \tau_2), \tau_2 \in T_0
\}$
are admissible. Hence, by Proposition \ref{P1.Adm} the families
$\{u_1(\theta), \theta\in T_S\}$ and $\{u_2(\theta),\theta\in
T_S\}$ are
admissible. Put
%
\begin{equation}\label{phi}
\phi(\theta)=\max[u_1(\theta), u_2(\theta)].
\end{equation}
The family $\{\phi(\theta), \theta\!\!\in T_S\}$, which is called
the \textit{new
reward family}, is also~clear\-ly admissible. \noindent Consider the
value function associated with the \mbox{new reward}
%
\begin{equation}\label{uS}
u(S)=\operatorname{ess}\sup_{\theta\in T_S} E[\phi(\theta) |
\mathcal{F}_S]
\qquad\mbox{a.s.}
\end{equation}

\begin{Theorem}[(Reduction)]\label{T3}
Suppose that $\{\psi(\theta,S), \theta,S\in T_0\}$ is a
biadmissible family.
For each stopping time $S$, consider $v(S)$ defined by (\ref{vS}) and
$u(S)$ defined by (\ref{u12}), (\ref{phi}), (\ref{uS}). Then,
\[
v(S)=u(S)\qquad \mbox{ a.s.}
\]
\end{Theorem}

\begin{pf}
Let $S$ be a stopping time.

\textit{Step 1}: First, let us show that $
v(S)\leq u(S)  \mbox{ a.s.}
$

Let $\tau_1,\tau_2\in
T_S$. Put $A=\{\tau_1\leq\tau_2\}$. As $A$ is in $\mathcal
{F}_{\tau_1}
\cap
\mathcal{F}_{\tau_2}$, we have
\[
E[\psi(\tau_1,\tau_2) |\mathcal{F}_S]= E[ {\mathbf{1}}_AE[\psi(\tau_1,\tau_2) |\mathcal{F}_{\tau_1}] |\mathcal
{F}_S] +E[{\mathbf{1}}_{A^
c}E[\psi(\tau_1,\tau_2) |\mathcal{F}_{\tau_2}] |\mathcal{F}_S].
\]
By noticing that on $A$ we have $E[\psi(\tau_1,\tau_2)
|\mathcal{F}_{\tau_1}] \leq u_2({\tau_1})\leq\phi({\tau_1\wedge
\tau
_2})$ a.s.~and, similarly, on $A^ c$ we have $E[\psi(\tau_1,\tau_2)
|\mathcal{F}_{\tau_2}] \leq u_1({\tau_2})\leq\phi({\tau_1\wedge
\tau
_2}) $ a.s., we~get
\[
E[\psi(\tau_1,\tau_2) |\mathcal{F}_S] \leq E[ \phi({\tau
_1\wedge\tau
_2})
|\mathcal{F}_S]\leq u(S)\qquad\mbox{a.s.}
\]
By taking the supremum over $\tau_1$ and $\tau_2$ in $T_S,$
we complete Step 1.

\textit{Step 2}: Let us now show that $v(S) \geq u(S)$ a.s.

We clearly have
$ v(S) \geq\operatorname{ess}\sup_{\tau_2 \in T_S}
E[\psi(S, \tau
_2 )
 |\mathcal{F}_S] = u_2 (S)$  a.s. By similar arguments, $v(S)
\geq
u_1(S)$ a.s.~and, consequently,
\[
v(S) \geq\max[ u_1(S), u_2 (S)] = \phi(S)\qquad \mbox{a.s.}
\]
Thus, $\{ v(S) ,S\in T_0\}$ is a supermartingale system which is
greater than $\{ \phi(S) , S\in T_0\}$. Now, by Proposition
\ref{P1.SuperM}, $\{ u(S) , S\in T_0\}$ is the smallest supermartingale
system which is greater than $\{ \phi(S) , S\in T_0\}$. Consequently,
Step 2 follows, which completes the proof.
\end{pf}

Note that the reduction to an optimal single stopping time
problem associated with a new reward will be the key property used to
construct optimal multiple stopping times and to establish an existence
result for them (see Sections \ref{2.3}--\ref{2.5}).

\subsection{Properties of optimal stopping times}\label{2.3}

$\!\!\!$In this section, we are given~a~biadmissible family
$\{\psi(\theta,S),\theta,S\!\!\in T_0\}$ such that \mbox{$ E[\operatorname{ess}
\sup_{\theta,S\in T_0} \psi(\theta,S)]\!<\!\infty$}.

\begin{Proposition}[(A necessary condition of optimality)]\label{P.CNopt}
Let $S$ be a~stopping time and consider the value function $v(S)$
defined by (\ref{vS}) for all $\theta\in T_S$, $u_1(\theta),
u_2(\theta)$ defined by (\ref{u12}), $\phi(\theta)$ defined by
(\ref{phi}) and $u(S)$ defined by (\ref{uS}).

Suppose that the pair $(\tau^*_1,\tau^*_2)$ is optimal for $v(S)$ and
put $A=\{\tau^*_1\leq\tau^*_2\}.$ Then:
\begin{enumerate}
\item[(1)]$\tau^*_1\wedge\tau^*_2$ is optimal for $u(S)$;
\item[(2)]
$\tau^*_2$ is optimal for $u_2({\tau^*_1})$ a.s.~on $A$;
\item[(3)]
$\tau^*_1$ is optimal for
$u_1({\tau^*_2})$ a.s.~on $A^c$.
\end{enumerate}
Moreover $A=\{\tau^*_1\leq\tau^*_2\}\subset B=\{
u_1(\tau^*_1\wedge\tau^*_2)\leq u_2(\tau^*_1\wedge\tau^*_2) \}$.
\end{Proposition}

\begin{pf}
Let $S$ $\in$ $T_0$ and suppose that the pair of
stopping times $(\tau^*_1, \tau^*_2)$ is optimal for $v(S)$. As
$u(S)=v(S)$ a.s., we obtain equality in Step 1 of the proof of Theorem
\ref{T3}. More precisely,
\begin{eqnarray*}
v(S)
&=&
E[\psi(\tau^*_1, \tau^*_2) |\mathcal{F}_S] =E[\phi(\tau^*_1\wedge\tau^*_2) |\mathcal{F}_S]=u(S)\qquad\mbox{a.s.},
\\
E[\psi(\tau^*_1,\tau^*_2)|\mathcal{F}_{\tau^*_1}]
&=&
u_2(\tau^*_1)=u_2(\tau^*_1\wedge\tau^*_2)=\phi(\tau^*_1\wedge\tau^*_2)\qquad\mbox{a.s. on } A,
\\
E[\psi(\tau^*_1,\tau^*_2)|\mathcal{F}_{\tau^*_2}]
&=&
u_1(\tau^*_2)=u_1(\tau^*_1\wedge\tau^*_2)=\phi(\tau^*_1\wedge\tau^*_2)\qquad\mbox{a.s.~on } A^c,
\end{eqnarray*}
which easily leads to (1), (2), (3) and $A \subset B$.
\end{pf}

\begin{Remark}\label{R2.1} Note that, in general, for a pair $(\tau
^*_1,\tau^*_2)$ of optimal stopping times for $v(S),$ the inclusion
$A\subset B$ is strict. Indeed if $\psi\equiv0$, then
$v=u=u_1=u_2=\phi=0$, and all pairs of stopping times are optimal.
Consider $\tau^*_1=T$, $\tau^*_2=0$. In this case, $A=\varnothing$ and
$B=\Omega$.
\end{Remark}

We now give a sufficient condition for optimality.

\begin{Proposition}[(Construction of optimal stopping times)]\label{Pconstruction}
Using the notation of Proposition \ref{P.CNopt},
suppose that:
\begin{enumerate}
\item$\theta^ *$ is optimal for $u(S)$;
\item$\theta^*_2$ is optimal
for $u_2({\theta^*})$;
\item
$\theta^*_1$ is optimal for
$u_1({\theta^*})$
\end{enumerate}
and put $B=\{ u_1(\theta^*) \leq u_2(\theta^ *)
\}$. The pair of stopping times $(\tau_1^*,\tau_2^*)$ defined~by
%
\begin{equation}\label{tau12}
\tau_1^*=\theta^{*}{\mathbf{1}}_{B}+ \theta^{*}_1{\mathbf{1}}_{B^c},
\qquad
\tau_2^*=\theta^*_2{\mathbf{1}}_{B}+ \theta^{*}{\mathbf{1}}_{B^c}
\end{equation}
is then optimal for $v(S)$.

Moreover, $\tau_1^*\wedge\tau_2^*=\theta^ *$ and $B= \{\tau
_1^*\leq
\tau_2^*\}.$
\end{Proposition}


\begin{pf}
Let $\theta^{*}$ be an optimal stopping time for
$u(S)$, that is, $u(S)=E[\phi({\theta^{*}}) |\mathcal
{F}_{S}]$ a.s. Let
$\theta^*_1$ be an optimal stopping time for $u_1({\theta^*})$ (i.e.,
$u_1({\theta^*})=E[\psi(\theta^*_1,\theta^{*}) |\mathcal
{F}_{\theta^*}]$
a.s.) and let $\theta^*_2$ be an optimal stopping time for
$u_2({\theta^*})$ (i.e., $u_2({\theta^*})=$
$E[\psi(\theta^{*}, \theta^*_2) |\mathcal{F}_{\theta^*}]$ a.s.).
We introduce the set $B=\{ u_1(\theta^*) \leq u_2(\theta^ *)
\}$. Note that $B$ is in $\mathcal{F}_{\theta^*}$.

Let $\tau_1^*,\tau_2^*$ be the stopping times defined by (\ref{tau12}).
We clearly have the inclusion
%
\begin{equation}\label{incl}
B \subset\{\tau_1^*\leq\tau_2^*\}.
\end{equation}
Since $u({S}) = E[\phi({\theta^{*}}) |\mathcal
{F}_{S}]$ and
$\phi({\theta^{*}}) = \max[u_1(\theta^*), u_2(\theta^ *)]$, we have
\[
u({S}) = E[{\mathbf{1}}_B u_2({\theta^*}) + {\mathbf{1} }_{B^c}
u_1({\theta^*}) |\mathcal{F}_S].
\]
The optimality of $\theta^*_1$ and $\theta^*_2$ gives that a.s.
\begin{eqnarray*}
u({S})
&=&
E[{\mathbf{1}}_B \psi(\theta^{*} ,\theta^*_2)+ {\mathbf{1}}_{B^c}\psi(\theta^*_1 ,\theta^{*}) |\mathcal{F}_S]
\\
&=&
E[{\mathbf{1}}_B \psi(\tau_1^*,\tau_2^*) + {\mathbf{1}}_{B^c}\psi(\tau_1^*,\tau_2^*) |\mathcal{F}_S]=E[\psi(\tau_1^*,\tau_2^*)|\mathcal{F}_S].
\end{eqnarray*}
As $u({S})= v(S)$ a.s., the pair of stopping times
$(\tau_1^*,\tau_2^*)$ is $S$-optimal for $v(S)$. By Proposition
\ref{P.CNopt}, we have $ \{\tau_1^*\leq\tau_2^*\} \subset B$.
Hence, by (\ref{incl}), $B = \{\tau_1^*\leq\tau_2^*\}$.
\end{pf}

\begin{Remark}
Proposition \ref{Pconstruction} still holds true if condition 2 holds
true on the set $B$ and condition 3 holds true on the set $B^c$.
\end{Remark}

Note that by Remark \ref{R2.1}, we do not have a characterization of
optimal pairs of stopping times. However, it is possible to give a
characterization of \textit{minimal optimal} stopping times in a
particular sense (see Appendix \hyperref[app2]{B}).

\subsection{Regularity of the new reward}\label{2.4}
Before studying the problem of the existence of optimal stopping
times, we have to state some regularity properties of the new reward
family $\{\phi(\theta),\theta\in T_0\}$.

Let us introduce the following definition.

\begin{Definition}
A biadmissible family $\{\psi(\theta,S),\theta,S\in T_0\}$ is said
to be \textit{uniformly right- (resp.,~left-) continuous in expectation along stopping
times [URCE (resp.,~ULCE)]} if $v(0)=\sup_{\theta,S\in T_0} E[\psi(\theta,S)]<\infty$ and
if, for each~$\theta$, $S\in T_0$ and
each sequence of stopping times $(S_n)_{n \in\mathbb{N}}$ such that
$S_n\downarrow S$ a.s.~(resp.,~$S_n\uparrow S$ a.s.),
\begin{eqnarray*}\label{hyp2P}
\lim_{n\to\infty}\sup_{ \theta\in T_0  } \vert E[   \psi (\theta, S)] - E[ \psi (\theta , S_n)  ] \vert &=& 0\quad \mbox{and}
\\
\lim_{n\to\infty}\sup_{\theta\in T_0  } \vert  E[   \psi (S, \theta)] - E[\psi (S_n,\theta)] \vert &=&0.
\end{eqnarray*}
\end{Definition}

The following right continuity property holds true for the new reward
family.

\begin{Theorem} \label{T2.NRcadgE}
Suppose that the biadmissible family $\{\psi(\theta,S),\theta,S\in
T_0\}$ is
URCE (resp., both URCE and ULCE).
The family $\{\phi(S), S\in T_0\}$ defined by (\ref{phi}) is then RCE (resp., both RCE and LCE).
\end{Theorem}

\begin{pf}
As $\phi(\theta)=\max[u_1(\theta),u_2(\theta)]$, it is sufficient
to show the RCE (resp., both RCE and LCE) properties for the family $\{u_1(\theta
),\theta\in T_0
\}$.\vadjust{\eject}

Let us introduce the following value function for each $S,\theta$ $\in$
$T_0$:
%
\begin{equation}\label{U1}
U_1 (\theta, S) = \operatorname{ess}\sup_{\tau_1 \in T_{\theta}}
E[\psi(\tau_1,
S)
|\mathcal{F}_{\theta}]\qquad\mbox{a.s.}
\end{equation}
As for all $\theta\in T_0$,
\[
u_1(\theta)=U_1(\theta,\theta)\qquad\mbox{a.s.},
\]
it is sufficient to prove that $\{U_1 (\theta, S), \theta, S \in
T_0\}$ is RCE (resp., both RCE and LCE), that is, if $\theta,S\in T_0$ and
$(\theta_n)_n$, $(S_n)_n$ in $T_0$ are such that $\theta_n\downarrow
\theta$ and $S_n\downarrow S$ a.s. (resp., $\theta_n\uparrow\theta$ and
$S_n\uparrow S$ a.s.), then $\lim_{n \to\infty} E[U_1
(\theta_n, S_n)] = E[U_1 (\theta, S)]$. Now, we have
\begin{eqnarray*}
&&
| E[U_1 (\theta, S)] -  E[U_1 (\theta_n, S_n)]
|
\\
&&\qquad\leq \underbrace{| E[U_1 (\theta, S)] - E[U_1 (\theta_n , S)]
|}_{\mbox{(I)}} + \underbrace{| E[U_1 (\theta_n, S)] - E[U_1 (\theta_n,
S_n)] |}_{\mbox{(II)}}.
\end{eqnarray*}
Let us show that (I) tends to $0$ as $n\to\infty$. For
each $S\in T_0$, $\{\psi(\theta,S),\theta\in T_0\}$ is an
admissible family
of positive random variables which is RCE (resp., both RCE and LCE). By Proposition
\ref{L1} (resp.,~Proposition \ref{P.vLCE}), the value function $\{U_1
(\theta, S), \theta\in T_0\}$ is RCE (resp., both RCE and LCE). It follows that (I)
converges to $0$ as $n$ tends to $\infty$.

Let us show that (II) tends to $0$ as $n\to\infty$. By
definition of the value function $U_1(\cdot,\cdot)$ (\ref{U1}), it follows that
\[
\vert E[U_1 (\theta_n, S)] - E[U_1 (\theta_n , S_n)] \vert
 \leq \sup_{\tau \in T_0 }  \vert E[   \psi (\tau, S)] - E[\psi (\tau , S_n)  ] \vert.
\]
which converges to $0$ since $\{\psi(\theta,S),\theta,S\in T_0\}
$ is URCE (resp., both URCE and ULCE). The proof of Theorem \ref{T2.NRcadgE} is thus complete.
\end{pf}

%

\begin{Corollary}
Suppose that
$v(0)= \sup_{\theta,S\in T_0}E[\psi(\theta,S)] < \infty$.
Under the same hypothesis as Theorem \ref{T2.NRcadgE}, the family
$\{v(S),S\in T_0\}$ defined by~(\ref{vS}) is RCE (resp., both RCE and LCE).
\end{Corollary}

\begin{pf}
This follows from the fact that $v(S)=u(S)$
a.s.~(Theorem \ref{T3}), where $\{u(S),S\in T_0\}$ is the value function
family associated with the new reward $\{\phi(S),S\in T_0\}$. Applying Propositions \ref{L1} and \ref{P.vLCE}, we obtain the required properties.
\end{pf}

We will now turn to the problem of the existence of optimal stopping~times.

\subsection{Existence of optimal stopping times}\label{2.5}


Let $\{\psi(\theta,S), \theta,S\in T_0\}$ be a~biadmissible
family which is
URCE and ULCE. Suppose that $v(0) < \infty$.\vadjust{\eject}

By Theorem \ref{T2.NRcadgE}, the admissible family of positive random
variables $\{\phi(\theta),\break \theta\!\in T_0\}$ defined by (\ref
{phi}) is RCE and
LCE. By Theorem \ref{T.1}, the stopping~time
\[
\theta^{*}=\operatorname{ess}\inf\{\theta\in T_S ,  u(\theta
)=\phi(\theta
) \mbox{ a.s.} \}
\]
is optimal for $u(S)\,[=v(S)]$, that is,
\[
u(S)=\operatorname{ess}\sup_{\theta\in T_S} E[\phi({\theta}) |
\mathcal{F}_{S}] =
E[\phi(\theta^{*}) |\mathcal{F}_{S}]\qquad\mbox{a.s.}
\]
Moreover, the
families $\{\psi(\theta,\theta^*), \theta\in T_{\theta^*}\}$ and
$\{\psi(\theta^*,\theta), \theta\in T_{\theta^*}\}$ are
admissible and are RCE and LCE. Consider the following optimal stopping
time problems defined for each $S \in T_{\theta^*}$:
\[
v_1(S)=\operatorname{ess}\sup_{\theta\in T_S}E[\psi(\theta,\theta^*)|
\mathcal{F}_{S}]
 \quad\mbox{and}\quad
v_2(S)=\operatorname{ess}\sup_{\theta\in T_S}E[\psi(\theta^*,\theta)|
\mathcal{F}_{S}].
\]
By
Theorem \ref{T.1} the stopping times $\theta^*_1$ and $\theta^*_2$
defined by $\theta^{*}_1=\operatorname{ess}\inf\{\theta\in T_{
\theta^{*}},\break
v_1(\theta)=\psi(\theta,\theta^*)  \mbox{ a.s.} \}$ and
$\theta^{*}_2=\operatorname{ess}\inf\{\theta\in T_{ \theta
^{*}} ,
v_2(\theta)=\psi(\theta^*,\theta)  \mbox{ a.s.} \}$ are optimal
stopping times for $v_1(\theta^*)$ and $v_2(\theta^*)$, respectively.
Note that $v_1(\theta^*)=u_1(\theta^*)$ and
$v_2(\theta^*)=u_2(\theta^*)$ a.s.


Let $\tau^*_1$ and $\tau^*_2$ be the stopping times defined by
%
\begin{equation}\label{mio}
\tau^*_1=\theta^*{\mathbf{1} }_B+ \theta^{*}_1{\mathbf{1} }_{B^c},
\qquad
\tau^*_2=\theta^*{\mathbf{1} }_{B^c}+ \theta^{*}_2{\mathbf{1} }_{B},
\end{equation}
where $B = \{ u^1(\theta^*)\leq u^2(\theta^*)\}$.
By Proposition \ref{Pconstruction}, the pair $(\tau_1^*, \tau_2^*)$
is optimal for $v(S)$. Consequently, we have proven the following theorem.

\begin{Theorem}[(Existence of an optimal pair of stopping times)]\label{Topt}
Let $\{\psi(\theta,S),\theta,S\in T_0\}$ be a biadmissible
family which is
URCE and ULCE. Suppose that $v(0) < \infty$.

The pair of stopping times $(\tau_1^*, \tau_2^*)$ defined by (\ref
{mio}) is then optimal for~$v(S)$ defined by (\ref{vS}).
\end{Theorem}

\begin{Remark}
Note that since $\theta^*$, $\theta^{*}_1$, $\theta^{*}_2$ are minimal
optimal, by results in Appendix \hyperref[app2]{B}, $(\tau_1^*, \tau_2^*)$ is minimal
optimal for $v(S)$ (in the sense defined in Appendix \hyperref[app2]{B}).
\end{Remark}


\section{The optimal $d$-stopping time problem}\label{3}
Let $d\in\mathbb{N}$, $d\geq2$. In this section, we show that
computing the
value function for the optimal $d$-stopping time problem
\[
v(S)=\operatorname{ess}\sup\{E[\psi( \tau_1,\ldots,\tau_d) |
\mathcal{F}_S], \tau
_1,\ldots,\tau_d\in T_S\}
\]
reduces to computing the value function for an optimal single stopping
time problem, that is,
\[
v(S)= \operatorname{ess}\sup\{E[\phi(\theta) |\mathcal
{F}_S], \theta\in T_S \}
\qquad\mbox{a.s.,}
\]
for a \textit{new reward} $\phi$. This new reward is expressed in terms
of optimal $(d-1)$-stopping time problems. Hence, by
induction, the initial optimal $d$-stopping time problem can
be reduced to \textit{nested} optimal single stopping time problems.

\subsection{Definition and initial properties of the value function}\label{3.1}

\begin{Definition}
We say that the family of random variables
$\{ \psi(\theta), \theta\in T_0^d\}$ is a $d$-\textit{admissible
family} if it
satisfies the following conditions:
\begin{enumerate}
\item for all $ \theta =  (\theta_1,  \ldots,\theta_d)\in
T_0^d,$ $\psi(\theta)$ is an $\mathcal{F}_{\theta_1  \vee
\cdots
\vee
\theta_d}$ measurable $\overline {\mathbb{R}}^+$-valued random variable;
\item
for all $ \theta, \theta'\in
T_0^d,\ \psi( \theta)=\psi( \theta')$ a.s.~on $
\{\theta=\theta'\} $.
\end{enumerate}
\end{Definition}

For each stopping time $S \in T_0$, we consider the value
function associated with the reward $\{\psi(\theta), \theta\in
T_0^d\}$:
%
\begin{equation}\label{vSd}
v(S)=\operatorname{ess}\sup_{\tau\in T_S^d} E[\psi(\tau) |
\mathcal{F}_S].
\end{equation}
As in the optimal double stopping time problem, the value
function satisfies the following properties.

\begin{Proposition}\label{P3.1} Let $\{\psi(\theta),\theta\in
T_0^d\}$ be
a $d$-admissible family of random variables. The following properties
then hold:
\begin{enumerate}
\item
$\{v(S),S\in T_0\}$ is an admissible family of random variables;
\item
For each $S\in T_0$, there exists a sequence of stopping times $(\theta
^n)_{n\in\mathbb{N}}$ in $T_S^d$
such that the sequence $\{ E[\psi(\theta^n) |\mathcal{F}_S]\}
_{ n\in
\mathbb{N}}$ is nondecreasing and such that $
v(S)=\lim_{n\to\infty} \uparrow E[\psi(\theta^{n}) |\mathcal
{F}_{S}]$ a.s.;
\item The family of random variables $\{ v(S) ,$ $ S\in T_0\}$ defined
by (\ref{vSd}) is a supermartingale system.
\end{enumerate}
\end{Proposition}

The proof is an easy generalization of the optimal double stopping time
problem (Proposition \ref{P2.1}).

Following the case with single or double stopping time, we now state the
following result on the regularity of the value function.

\begin{Proposition}\label{Ld.vRCE}
Suppose that the $d$-admissible family $\{\psi(\theta), \theta\in
T_0^d\}$
is RCE and that  $v(0)<\infty$.
The family $\{ v(S) ,
S\in T_0
\}$ defined by (\ref{vSd}) is then RCE.
\end{Proposition}

The definition of RCE and the proof of this property are
easily derived from the single or double stopping time case (see Definition
\ref{RCE2}
and Proposition \ref{L2.vRCE}).

\subsection{Reduction to an optimal single stopping time problem}\label{3.2}
The optimal $d$-stopping time problem (\ref{vSd}) can be
expressed in terms of an optimal single stopping time problem as
follows.

For $i=1,\ldots,d$ and $\theta\in T_0$, consider the random variable
%
\begin{eqnarray}\label{uiun}
&&\quad\quad u^{(i)}(\theta) =\operatorname{
ess}
\sup_{\tau_1,\ldots,\tau_{i-1},\tau_{i+1},\ldots
,\tau_{d}
\in T_{\theta}^{d-1}}
E[\psi(\tau_1,\ldots,\tau_{i-1},\nonumber
\\[-8pt]\\[-8pt]
&&\quad\quad\hphantom{u^{(i)}(\theta) =\operatorname{
ess}
\sup_{\tau_1,\ldots,\tau_{i-1},\tau_{i+1},\ldots
,\tau_{d}
\in T_{\theta}^{d-1}}
E[\psi(} \theta,\tau_{i+1},\ldots,\tau
_{d})
|\mathcal{F}_{\theta}].\nonumber
\end{eqnarray}
Note that this notation is adapted to the $d$-dimensional case.

In the two-dimensional case ($d=2$), we have
\[
u^{(1)}(\theta)= \operatorname{ess}\sup_{\tau_2\in
T_{\theta}} E[\psi(\theta,\tau_2) |\mathcal{F}_{\theta}]=
u_2({\theta})\qquad
\mbox{a.s.}
\]
and
\[
u^{(2)}(\theta)= \operatorname{ess}\sup_{\tau_1\in T_{\theta}}
E[\psi(\tau_1,
\theta) |\mathcal{F}_{\theta}] = u_1({\theta}) \qquad\mbox{a.s.},
\]
by definition of $u_1({\theta})$ and $u_2({\theta})$ [see
(\ref{u12})]. Thus, the notation in the two-dimensional case was
different, but more adapted to that simpler case.

For each $\theta\in T_0$, define the $\mathcal{F}_{\theta
}$-measurable random
variable called the \textit{new reward},
%
\begin{equation}\label{d-phi}
\phi(\theta)=\max\bigl[u^{(1)}(\theta),\ldots, u^{(d)}(\theta)\bigr],
\end{equation}
and for each stopping time $S$, define the $\mathcal{F}_S$-measurable variable
%
\begin{equation}\label{d-uS}
u(S)=\operatorname{ess}\sup_{\theta\in T_S} E[\phi(\theta) |
\mathcal{F}_S] .
\end{equation}

\begin{Theorem}[(Reduction)] \label{T3.red}
Let $\{\psi(\theta), \theta\in T_0^d\}$ be a $d$-admissible family
of random variables and for each stopping time $S$, consider $v(S)$
defined by~(\ref{vSd}) and $u(S)$ defined by~(\ref{uiun}),
(\ref{d-phi}) and (\ref{d-uS}). Then,
\[
v(S)=u(S)\qquad \mbox{a.s.}
\]
\end{Theorem}

\begin{pf}
\textit{Step 1}: Let us prove that for all $S\in T_0$, $
v(S)\leq
u(S)$ a.s.

Let $S$ be a stopping time and $\tau=(\tau_1,\ldots,\tau_d)\in T_S^d$.
There exists $(A_i)_{i=1,\ldots,d}$ with $\Omega= \bigcup_{i} A_i$, where
$A_i \cap A_j = \empty$ for $i \neq j$,
$\tau_1\wedge\cdots\wedge\tau_d=\tau_i$ a.s.~on $A_i$ and~$A_i$
are in
$\mathcal{F}_{\tau_1\wedge\cdots\wedge\tau_d}$ for $i=1,\ldots,
d$ (for $d=2$,
one can take $A_1=\{\tau_1\leq\tau_2\}$ and $A_2=A_1^c$). We
have
\[
E[\psi(\tau) |\mathcal{F}_S]=\sum_{i=1}^dE[ {\mathbf{1}}_{A_i}E[\psi(\tau) |\mathcal{F}_{\tau_i}] |\mathcal{F}_S].
\]
By noticing that
on $A_i$ one has a.s.~$E[\psi(\tau) |\mathcal{F}_{\tau_i}] \leq
u^{(i)}({\tau_i})\leq\phi(\tau_i)=
\phi(\tau_1\wedge\cdots\wedge \tau_d)$, we get
${E[\psi(\tau) |\mathcal{F}_S] \leq E[
\phi({\tau_1\wedge\cdots\wedge \tau_d})  |\mathcal{F}_S]\leq
u(S)
\mbox{ a.s.}}$ By taking the supremum over $\tau=
(\tau_1,\ldots,\tau_d),$ we complete Step 1.

\textit{Step 2}:  Let us show that for all $S\in T_0$,
$v(S)\geq u(S)$ a.s.

This follows from the fact that $\{v(S),S\in T_0\}$ is a
supermartingale system greater than $\{\phi(S),S\in T_0\}$ and
that $\{
u(S),S\in T_0\}$ is the smallest supermartingale system of this
class.\vspace*{-2pt}
\end{pf}

Note that the new reward is expressed in terms of optimal $(d-1)$-stopping time problems. Hence, by induction, the initial
optimal $d$-stopping time problem can be reduced to nested
optimal single stopping time problems. In the case of a symmetric reward,
the problem reduces to ordered stopping times and the nested optimal single stopping time problems simply reduce to a sequence of optimal
single
stopping time problems defined by backward induction (see Section
\ref{swings} and the application to \textit{swing options}).\vspace*{-2pt}

\subsection{Properties of optimal stopping times in the $d$-stopping
time problem}\label{3.3}

Let $\{\psi(\theta), \theta\in T_0^d\}$ be a $d$-admissible
family. Let us introduce the following notation: for
$i=1,\ldots,d$, $\theta\in T_0$ and $\tau_1,\ldots,\tau_{d-1}$ in
$T_0$, consider
the random variable
%
\begin{equation}\label{psii}
\psi^{(i)}(\tau_1,\ldots,\tau_{d-1},\theta)
=\psi(\tau_1,\ldots,\tau_{i-1},\theta,\tau_{i},\ldots,\tau_{d-1}).
\end{equation}
Using this notation, note that for each $i=1,\ldots,d,$ the value
function $u^{(i)}$ defined at (\ref{uiun}) can be written
%
\begin{equation}\label{ui}
u^{(i)}(\theta) =\operatorname{ess}\sup_{\tau\in T_{\theta}^{d-1}}
E\bigl[\psi^{(i)}(\tau,\theta) |\mathcal{F}_{\theta}\bigr].\vspace*{-2pt}
\end{equation}

\begin{Proposition}[(Construction of optimal stopping times)]\label{d-construction}
$\!\!$Suppose~that:\vspace*{-5pt}
\begin{enumerate}
\item there exists an optimal stopping time $\theta^ *$ for $u(S)$;
\item for $i=1,\ldots,d$, there exist
$(\theta^{(i)*}_1,\ldots,\theta^{(i)*}_{i-1},\theta
^{(i)*}_{i+1},\ldots,\theta^{(i)*}_d)=\theta^{(i)*}$
in $T_{\theta}^{d-1}$ such that
$u^{(i)}({\theta^*})=E[\psi^{(i)}(\theta^{(i)*},\theta^*)
|\mathcal{F}_{\theta^*}].$
\end{enumerate}
Let $(B_i)_{i=1,\ldots,d}$ with $\Omega= \bigcup_{i} B_i$ be such that
$B_i \cap B_j = \varnothing$ for $i \neq j$, $\phi(\theta^*)=$
$u^{(i)}(\theta^*)$ a.s.~on $B_i$ and $B_i$ is
$\mathcal{F}_{\theta^*}$-measurable for $i=1,\ldots,d$. Put\vspace*{-1pt}
%
\begin{equation}\label{tauj}
\tau_j^*=\theta^* {\mathbf{1} }_{B_j}+ \sum_{i\neq j, i=1}^d\theta
^{(i)*}_j{\mathbf{1}
}_{B_i}.
\end{equation}
Then, $(\tau^*_1,\ldots,\tau^*_d)$ is optimal for $v(S)$, and
$\tau_1^*\wedge\cdots\wedge\tau_d^*= \theta^*$.\vspace*{-2pt}
\end{Proposition}

\begin{pf}
It is clear that
$\tau_1^*\wedge\cdots\wedge\tau_d^*= \theta^*$, and a.s.\vspace*{-1pt}
\begin{eqnarray*}
u({S})
& =&
E[\phi(\theta^*) |\mathcal{F}_S]=\sum_{ i=1}^dE\bigl[{\mathbf{1}}_{B_i} u^{(i)}({\theta^*})|\mathcal{F}_S\bigr]
\\
&= &
\sum_{ i=1}^dE\bigl[{\mathbf{1}}_{B_i} E\bigl[ \psi^{(i)}\bigr(\theta^{(i)*},\theta^*\bigr)|\mathcal{F}_{\theta^*}\bigr]|\mathcal{F}_S\bigr]
\\
&=&
\sum_{ i=1}^dE\bigl[{\mathbf{1}}_{B_i}E\bigl[\psi\bigl(\theta^{(i)*}_1,\ldots,\theta^{(i)*}_{i-1},\theta^*,\theta^{(i)*}_{i+1},\ldots,\theta^{(i)*}_d\bigr)
|\mathcal{F}_{\theta^*}\bigr]|\mathcal{F}_S\bigr]
\\
&=&
E[\psi(\tau^*_1,\ldots,\tau^*_{i-1},\tau^*_{i},\tau
^*_{i+1},\ldots,\tau^*_d)
|\mathcal{F}_S] \leq v(S)=u(S).
\end{eqnarray*}
\upqed\end{pf}

\begin{Remark}
As in the bidimensional case, one can easily derive a~necessary
condition for obtaining optimal stopping times. Moreover, for an
adapted partial order relation on $\mathbb{R}^d$, one can also derive a
characterization of minimal optimal $d$-stopping times. This result is
given in Appendix~\hyperref[app22]{B.2}.\vspace*{-1pt}
\end{Remark}

Before studying the existence of an optimal $d$-stopping time for
$v(S)$, we will study the regularity properties of the new reward
$\{\phi(\theta),\theta\in T_0\}$ defined by (\ref{d-phi}).\vspace*{-1pt}

\subsection{Regularity of the new reward}\label{3.4}

Let us introduce the following definition of uniform continuity.\vspace*{-1pt}

\begin{Definition} A
$d$-admissible family $\{\psi(\theta),\theta\in T_0^d\}$ is said
to be \textit{uniformly right- (resp.,~left-) continuous along stopping times in
expectation [URCE (resp.,~ULCE)]} if $
v(0)<\infty$, and for each
$i=1,\ldots,d$, $S\in T_0$ and sequence of stopping times $(S_n)_{n \in
\mathbb{N}}$
such that $S_n\downarrow S$ a.s.~(resp.,~$S_n\uparrow S$ a.s.), we have\vspace*{-1pt}
\[
\lim_{n\to\infty} \sup_{  \theta\in  T_0^{d-1} }\bigl| E\bigl[   \psi^{(i)} (\theta, S_n)\bigr] - E\bigl[\psi^{(i)} (\theta, S)\bigr]\bigr|  = 0 \qquad
\mbox{a.s.}\vspace*{-1pt}
\]
\end{Definition}

\begin{Proposition} \label{Pd.NGRCE}
Let $\{\psi(\theta),\theta\in T_0^d\}$ be a $d$-admissible
family which is
URCE (resp., both URCE and ULCE). The family of positive random variables $\{\phi(S),
S\in T_0\}$ defined by (\ref{d-phi}) is then RCE (resp., both RCE and LCE).\vspace*{-1pt}
\end{Proposition}

\begin{pf}
The proof uses an induction argument. For $d=1$
and $d=2$, the result has already been shown. Fix $d\geq1$ and
suppose by induction that the property holds for any $d$-admissible
family which is URCE (resp., both URCE and ULCE). Let $\{\psi(\theta),\theta\in
T_S^{d+1}\}$ be a $(d+1)$-admissible family which is URCE (resp., both URCE and ULCE).
As $\phi(\theta)=\max[u^{(1)}(\theta),\ldots,\break u^{(d+1)}(\theta)]$,
it is sufficient to show the RCE (resp., both RCE and LCE) properties for the family $\{
u^{(i)}(\theta),\theta\in T_0\}$ for all $i=1,\ldots,d+1$.

Let us introduce the following value function for each $S,\theta$ $\in$
$T_0$:
%
\begin{equation}\label{Ui}
U^{(i)} (\theta, S) = \operatorname{ess}\sup_{\tau\in T_{\theta
}^{d}} E\bigl[\psi^{(i)}
(\tau, S) |\mathcal{F}_{\theta}\bigr]\qquad\mbox{a.s.}
\end{equation}
As for all $\theta\in T_0$,\vspace*{-1pt}
\[
u^{(i)}(\theta)=U^{(i)}(\theta,\theta)
\qquad\mbox{a.s.},
\]
it is sufficient to prove that the biadmissible family $\{
U^{(i)}(\theta, S),\theta,S\in T_0\}$ is RCE (resp., both RCE and LCE) as in the
bidimensional case.\vadjust{\eject}

Let $\theta,S\in T_0$ and $(\theta_n)_n$, $(S_n)_n$ be monotonic
sequences of stopping times that converge, respectively, to $\theta$
and $S$ a.s. We have
\begin{eqnarray*}
&&E\bigl[\bigl| U^{(i)} (\theta, S)-  U^{(i)} (\theta_n, S_n) \bigr|\bigr]
\\
&&\qquad\leq
  \underbrace{E\bigl[\bigl| U^{(i)} (\theta, S) - U^{(i)} (\theta_n , S)
\bigr|\bigr]}_{\mbox{(I)}} + \underbrace{E\bigl[\bigl| U^{(i)} (\theta_n, S) - U^{(i)}
(\theta_n, S_n) \bigr|\bigr]}_{\mbox{(II)}}.
\end{eqnarray*}
Let us show that (I) tends to $0$ as $n\to\infty$. Note that for each
$S\in T_0$, $\{\psi^{(i)}(\tau,S),\tau\in T_0^{d}\}$ is a
$d$-admissible family of positive random variables which is URCE (resp., both URCE and ULCE) and $\{U^{(i)} (\theta, S)$, $\theta\in T_0\}$ is the corresponding
value function family. By the induction assumption, this family is RCE
(resp., both RCE and LCE). Hence, (I) converges a.s.~to $0$ as~$n$ tends to $\infty$ when
$(\theta_n)$ is monotonic.

Let us now show that (II) tends to $0$ as $n\to\infty$. By
definition of the value function $U^{(i)}(\cdot,\cdot)$ (\ref{Ui}), it follows
that
\[
E\bigl[\bigl| U^{(i)} (\theta_n, S) - U^{(i)} (\theta_n , S_n) \bigr|\bigr]
\leq  \sup_{\theta \in T_0^{d}  } \bigl| E\bigl[ \psi^{(i)} (\theta, S)\bigr] - E\bigl[\psi^{(i)} (\theta , S_n)\bigr]
\bigr|,
\]
and the right-hand side tends to $0$ by the URCE (resp., both URCE and ULCE) properties of
$\psi$.
\end{pf}

\subsection{Existence of optimal stopping times}\label{3.5}

By Theorem \ref{T.1}, the regularity properties of the new reward will
ensure the existence of an optimal stopping time $\theta^*\in T_0$ for
$u(S)$. By Proposition \ref{d-construction}, this will allow us to
show by induction the existence of an optimal stopping time for $v(S)$.

\begin{Theorem}[(Existence of optimal stopping times)]\label{Toptd}
Let $\{\psi(\theta), \theta\in T_0^d\}$ be a
$d$-admissible family
of positive random variables which is URCE and ULCE. There then exists
a $\tau^*\in T_S^d$ optimal for $v(S)$, that is, such that
\[
v(S)=\operatorname{ess}\sup_{\tau\in T_S^d} E[\psi(\tau)
|\mathcal{F}_S]= E[\psi(\tau
*)|\mathcal{F}_S].
\]
\end{Theorem}

\begin{pf}
The result is proved by
induction on $d$.
For $d=1$ the result is just Theorem \ref{T.1}. Suppose now that
$d\geq1$ and suppose by induction that for all $d$-admissible
families which are URCE and ULCE, optimal $d$-stopping times do
exist.
Let $\{\psi(\theta), \theta\in T_S^{d+1}\}$ be a
$(d+1)$-admissible family which is URCE and ULCE. The
existence of an optimal $(d+1)$-stopping time for the associated value
function $v(S)$ will be derived
by applying Proposition~\ref{d-construction}.
Now, by Proposition~\ref{Pd.NGRCE}, the new reward family $\{\phi
(\theta), \theta\in T_0\}$ is
LCE and RCE. By Theorem~\ref{T.1}, there exists an optimal stopping
time $\theta^*$ for $u(S)$. Thus, we have proven that condition 1 of
Proposition \ref{d-construction} is satisfied.

Note now that for $i=1,\ldots,d+1$, the $d$-admissible families $\{
\psi^{(i)}(\theta,\theta^*),\theta\in T_0^d\}$ are URCE and ULCE.
Thus, by the
induction hypothesis, for each $\theta\in T_0$, there exists an
optimal $\theta
^{*(i)}\in T_{\theta^*}^d$ for the value function\vadjust{\eject}
$U^{(i)}(\theta^*,\theta^*) $ defined by (\ref{Ui}).
Noting that $U^{(i)}(\theta^*,\theta^*)=u^{(i)}(\theta^*)$, we have
proven that condition 2 of Proposition \ref{d-construction} is
satisfied. Now applying Proposition \ref{d-construction}, the result follows.
\end{pf}

\subsection{Symmetric case}\label{swings}

Suppose that $\psi(\tau_1,  \ldots,\tau_d)$ is symmetric with
respect to $(\tau_1,  \ldots,\tau_d),$ that is,
\[
\psi(\tau_1,  \ldots,\tau_d) = \psi\bigl(\tau_{\sigma(1)},
\ldots,\tau_{\sigma(d)}\bigr)
\]
for each permutation $\sigma$ of $\{1,\ldots, d\}$. By symmetry we can
suppose that $\tau_1 \leq\tau_2\leq\cdots\leq\tau_d$, that is, that
the value function $v(S)$ coincides with
\[
{ v_d(S)=\operatorname{ess} \sup_ { (\tau
_1,\ldots,\tau_d) \in\mathcal{S}^d_S } E[\psi( \tau_1,\ldots,\tau
_d) |\mathcal{F}_S]},
\]
where $\mathcal{S}^d_S = \{ \tau_1,\ldots, \tau_d \in T_S$ s.t. $\tau_1 \leq\tau_2\leq\cdots\leq\tau_d \}$.
It follows that the value functions $u^{(i)}(\theta) $ and the new
reward $\phi(\theta)$ coincide and are simply given for each $\theta
\in T_0$ by
the following random variable:
\[
\phi_1(\theta) =\operatorname{ess} \sup_ {
(\tau_2 , \tau_3, \ldots, \tau_{d}) \in\mathcal{S}^{d-1}_{\theta}}
E[\psi(\theta, \tau_2,\ldots, \tau_{d}) |\mathcal{F}_{\theta}].
\]
The reduction property can be written as follows:
\[
v(S)=\operatorname{ess}\sup_{\theta\in T_S} E[\phi_1(\theta) |
\mathcal{F}_S].
\]
We then consider the value function $\phi_1(\theta_1)$.
The associated new reward is given for $\theta_1$, $\theta_2$ such that
$S\leq\theta_1 \leq\theta_2$ by
\[
\phi_2(\theta_1, \theta_2) =\operatorname{ess} \sup_ {
 (\tau_3,\ldots, \tau_{d}) \in\mathcal{S}^{d-2}_{\theta_2}}
E[\psi(\theta_1, \theta_2, \tau_3,\ldots, \tau_{d})
|\mathcal{F}_{\theta_2}].
\]
Again, the reduction property gives
%
\begin{equation}\label{phiun}
\phi_1(\theta_1)=\operatorname{ess}\sup_{\theta\in T_{ \theta_1 }
} E[\phi
_2(\theta_1, \theta_2) |\mathcal{F}_{\theta_1}].
\end{equation}
We then consider the value function $\phi_2(\theta_1, \theta_2),$ and
so on. Thus, by forward induction, we define the new rewards $\phi_{i}$
for $i=1,2, \ldots, d-1$ by
\[
\phi_{i}( \theta_1,\ldots, \theta_{i})=\operatorname{ess} \sup_
{ (\tau_{i+1},\ldots, \tau_{d}) \in\mathcal{S}^{d-i}_{\theta_i}} E[\psi(\theta_1,\ldots, \theta_i,
\tau_{i+1},\ldots, \tau_{d}) |\mathcal{F}_{\theta_i}]
\]
for each $(\theta_{1},\ldots, \theta_{i}) \in\mathcal{S}^{i}_{S}$. The
reduction property gives
%
\begin{equation}\label{phii}
\phi_{i}( \theta_1,\ldots, \theta_{i})=\operatorname{ess}\sup_{
\theta_{i+1} \in T_{
\theta_{i}} } E[\phi_{i+1}( \theta_1,\ldots, \theta_{i},
\theta_{i+1}) |\mathcal{F}_{ \theta_{i} } ].
\end{equation}
Note that for $i=d-1$,
%
\begin{equation}\label{phid}
\phi_{d-1}( \theta_1,\ldots, \theta_{d-1})=\operatorname{ess}\sup
_{ \theta_{d}
\in
T_{ \theta_{d-1}} } E[\Psi( \theta_1,\ldots, \theta_{d-1},
\theta_{d}) |\mathcal{F}_{ \theta_{d-1} } ]
\end{equation}
for each $(\theta_{1},\ldots, \theta_{d-1})
\in\mathcal{S}^{d-1}_{S}$.\vadjust{\eject}

Hence, using backward induction we can now define $\phi_{d-1}(
\theta_1,\ldots, \theta_{d-1})$ by (\ref{phid}) and then $\phi_{d-2}(
\theta_1,\ldots, \theta_{d-2}),\ldots,\phi_2(\theta_1,
\theta_2),\phi_1(\theta_1)$ by the induction formula (\ref{phii}).
Consequently, we have the following characterization of the value
function and construction of a multiple optimal stopping time (which
are rather intuitive). \vadjust{\goodbreak}

\begin{Proposition}\label{popo}
\begin{itemize}
\item
Let $\{\psi(\theta), \theta\in T_0^d\}$ be a symmetric
$d$-admissible family of random variables, and for each stopping time
$S$, consider the associated value function~$v(S)$.

Let $\phi_i$, $i=d-1,  d-2,\ldots, 2, 1,$ be defined by
backward\vspace*{1pt} induction as follows:
$\phi_{d-1}( \theta_1,\ldots, \theta_{d-1})$ is given by (\ref{phid})
for each $(\theta_{1},\ldots, \theta_{d-1}) \in\mathcal{S}^{d-1}_{S}$.
Also, for $i= d-2,\ldots,2,1$ and each $(\theta_{1},\ldots,
\theta_{i}) \in\mathcal{S}^{i}_{S}$, $\phi_{i}( \theta_1,\ldots,
\theta_{i})$ is given
in terms of the function $\phi_{i+1}$ by backward induction formula (\ref{phii}).

The value function then satisfies
%
\begin{equation}
\label{vf} v(S)=\operatorname{ess}\sup_{\theta\in T_S} E[\phi
_1(\theta) |\mathcal{F}_S].
\end{equation}
\item Suppose that $\{\psi(\theta), \theta\in T_0^d\}$ is URCE
and ULCE. Let
$\theta_1^*$ be an optimal stopping time for $v(S)$ given by
(\ref{vf}), let $\theta_2^*$ be an optimal stopping time for
$\phi_1(\theta_1^*)$ given by (\ref{phiun}) and for $i=2,3,\ldots, d-1$, let $\theta_{i+1}^*$ be an optimal stopping time for
$\phi_{i}(\theta_1^*,\ldots, \theta_{i}^*)$ given by (\ref
{phii}).

Then, $(\theta_1^*,\ldots, \theta_{d}^*)$ is a multiple optimal
stopping time for $v(S)$.
\end{itemize}
\end{Proposition}

\subsubsection*{Some simple examples}

First, consider the very simple additive case: suppose that the reward
is given by
%
\begin{equation}\label{add}
\psi (\tau_1, \ldots ,\tau_d)= Y({\tau_1}) + Y({\tau_2})+ \cdots + Y({\tau_d}),
\end{equation}
where $Y$ is an admissible family of random variables such that\break
 $\sup_{\tau \in T_0} E[Y(\tau)]\!<\! \infty$. We then obviously
have that
$v(S)\!=\! dv^1(S)$, where~$v^1(S)$ is the value function of the single
optimal stopping time problem associated with reward $Y$. Also, if
$\theta_1^*$ is an optimal stopping time for $v_1(S)$, then $(\theta
_1^*,\ldots, \theta_1^*)$ is optimal for $v(S)$.

\subsubsection*{Application to swing options}
Let us now consider the more interesting additive case of \textit{swing
options}: suppose that $T = + \infty$ and that the reward is still
given by (\ref{add}), but the stopping times are separated by a fixed
amount of time $\delta>0$ (sometimes called ``refracting time''). In
this case, the value function is given by
\[
v(S)=\operatorname{ess}\sup\{E[\psi( \tau_1,\ldots,\tau_d) |
\mathcal{F}_S], (\tau
_1,\ldots,\tau_d) \in\mathcal{S}^d_S\},
\]
where $\mathcal{S}^d_S = \{ \tau_1,\ldots, \tau_d \in T_S \mbox{ s.t. } \tau_i \in T_{\tau_{i-1} + \delta},2 \leq i \leq
d-1\}$.
All the previous properties then still hold. Again, the $\phi_i$
satisfy the following induction\vadjust{\eject} equality:\vspace*{-1pt}
\[
\phi_{i}( \theta_1,\ldots, \theta_{i})=\operatorname{ess}\sup_{
\theta_{i+1} \in T_{
\theta_{i}+ \delta} } E[\phi_{i+1}( \theta_1,\ldots, \theta_{i},
\theta_{i+1}) |\mathcal{F}_{ \theta_{i} } ].
\]
%
One can then easily derive that
$\phi_{d-1}( \theta_1,\theta_2,\ldots,  \theta_{d-1})=
 Y({ \theta_1}) + \cdots +Y({\theta_{d-1}})+  Z_{d-1}(  \theta_{d-1} )$, where\vspace*{-1pt}
\begin{eqnarray*}
Z_{d-1}( \theta_{d-1} )&=&\operatorname{ess}\sup_{ \tau\in T_{ \theta
_{d-1}+ \delta
} } E[Y(\tau) |\mathcal{F}_{ \theta_{d-1} } ].
\\[-1pt]
\phi_{d-2}( \theta_1,\ldots, \theta_{d-2})&=& Y(\theta_1) +
\cdots
+Y(\theta_{d-2})+ Z_{d-2}( \theta_{d-2} ),\qquad \mbox{where}
\\[-1pt]
Z_{d-2}( \theta_{d-2} )&=&\operatorname{ess}\sup_{ \tau\in T_{ \theta
_{d-2}+ \delta
} } E[Y(\tau)+Z_{d-1}( \tau)  |\mathcal{F}_{ \theta_{d-2} }
],
\end{eqnarray*}
and so on. Hence, for $i=1,2,\ldots, d-2$,
$\phi_{i}( \theta_1,\ldots, \theta_{i})=Y(\theta_1) + \cdots
+Y(\theta_{i})+ Z_{i}( \theta_{i} )$,
where\vspace*{-1pt}
\[
Z_{i}( \theta_{i} )=\operatorname{ess}\sup_{ \tau\in T_{ \theta
_{i}+ \delta} }
E[Y(\tau)+Z_{i+1}( \tau)  |\mathcal{F}_{ \theta_{i} } ].
\]
The value function satisfies\vspace*{-1pt}
%
\begin{equation}\label{vsdeu}
v(S)=\operatorname{ess}\sup_{\theta\in T_S} E[Y(\theta)+ Z_1(\theta)|\mathcal{F}_S].
\end{equation}
This corresponds to Proposition 3.2 of Carmona and Dayanik (\citeyear{CD}).

Suppose that $Y$ is RCE and LCE. Let $\theta_1^*$ be the minimal
optimal stopping time for $v(S)$ given by (\ref{vsdeu}) and for $i=1,2
,\ldots, d-1$, let $\theta_{i+1}^*$ be the minimal optimal stopping
time for $Z_{i}( \theta_{i}^*)$. The $d$-stopping time
$(\theta_1^*,\ldots, \theta_{d}^*)$ is then the minimal optimal
stopping time for $v(S)$. This corresponds to Proposition 5.4 of
Carmona and Dayanik (\citeyear{CD}).

Note that the multiplicative case can be solved similarly. Further
applications to American options with multiple exercise times are
studied in Kobylanski and Quenez (\citeyear{KQ}).\vspace*{-1pt}

\section{Aggregation and multiple optimal stopping times}\label{3.9}

As explained in the \hyperref[in]{Introduction}, in previous works on
the optimal single stopping time problem, the~reward is given by an
RCLL positive adapted process $(\phi_t) $. Moreover, when the reward
$(\phi_t)$ is continuous, an optimal $S$-stopping time is given by\vspace*{-1pt}
%
\begin{equation}\label{overtheta}
\overline\theta(S)=\inf\{t\geq S, v_t=\phi_t\},
\end{equation}
which corresponds to the first hitting time after $S$ of $0$ by the
RCLL adapted process $(v_t -\phi_t)$.
This formulation is very important since it gives a simple and
efficient method to compute an optimal stopping time.

In the two-dimensional case, instead of considering a reward process,
it is quite natural to suppose that the reward is given by a biprocess
$(\Psi_{t,s})_{(t,s) \in[0,T] ^2}$ such that a.s., the map $(t,s)
\mapsto\Psi_{t,s}$ is continuous and for each $(t,s) \in[0,T] ^2$,
$\Psi_{t,s}$ is $\mathcal{F}_{t \vee s}$-measurable (see Remark
\ref{biprocess}).\vadjust{\eject}

We would like to construct some optimal stopping times by using hitting
times of processes. By the existence and construction properties of
optimal stopping times given in Theorem \ref{Topt}, we are led to
construct $\theta^*$, $\theta^*_1$ and $\theta^*_2$ as hitting times
of processes. Since $\Psi$ is a continuous biprocess, there is no
problem for $\theta^*_1$, $\theta^*_2$. However, for $\theta^*$ we need
to aggregate the new reward $\{\phi(\theta),\theta\in T_0\}$,
which requires
new aggregation results. These results hold under stronger assumptions
on the reward than those made in the previous existence theorem
(Theorem~\ref{Topt}).

\subsection{Some general aggregation results}\label{3.10}
\subsubsection{Aggregation of a supermartingale system}\label{3.10.1}
Recall the classical result of aggregation of a supermartingale system
[El Karoui (\citeyear{EK})].

\begin{Proposition}\label{P.SMA}
Let $\{h(S),S\in T_0\}$ be a supermartingale system which is RCE and
such that $h(0) < \infty$. There then exists an RCLL adapted process
$(h_t)$ which \textit{aggregates} the family $\{h(S), S \in T_{0}\},$
that is, for each $S\in T_0$, $h_{S}= h(S) $ a.s.
\end{Proposition}

This lemma relies on a well-known result [see, e.g., El Karoui (\citeyear{EK})
or Theorem 3.13 in Karatzas and Shreve (\citeyear{KS1}); for details, see the
proof in Section \ref{3.14}].

Classically, the above Proposition \ref{P.SMA} is used to aggregate the
value function of the single stopping time problem. However, it cannot
be applied to the new reward since it is no longer a supermartingale
system. Thus, we will now state a new result on aggregation.


\subsubsection{A new result on aggregation of an admissible family}\label{3.11}

Let us introduce the following right-continuous property for admissible
families.

\begin{Definition}
An admissible family $\{\phi(\theta),\theta\in T_0\}$ is said to
be \textit{right-continuous along stopping times (RC)} if for any $\theta\in
T_0$ and any sequence $(\theta_n)_{n\in\mathbb{N}}$ of stopping
times such
that $\theta_n\downarrow\theta$ a.s., we have $
{\phi(\theta)=\lim_{n\to\infty} \phi(\theta_n)} $ a.s.
\end{Definition}

We state the following result.

\begin{Theorem} \label{T.Aggregphi}
Suppose that the admissible family of positive random variables $\{
\phi(\theta), \theta\in T_0 \}$ is right-continuous along stopping
times. There then exists a progressive process $( \phi_t ) $ such that
for each $\theta \in T_0$, $ \phi_{\theta}= \phi(\theta)$ a.s. and such that there exists a nonincreasing sequence
of right-continuous processes~$(\phi_t ^{n})_{n\in\mathbb{N}}$ such
that for
each $(\omega, t) \in \Omega\times[0,T]$, $\lim_{n \rightarrow
\infty} \phi_t ^{n}(\omega) =\break \phi_t (\omega)$.
\end{Theorem}

\begin{pf}
See Section \ref{3.14}.
\end{pf}

\subsection{The optimal stopping problem}\label{3.12}

First, recall the following classical result [El Karoui (\citeyear{EK})].\vspace*{-1pt}
\begin{Proposition}[(Aggregation of the value function)]\label{P1.Aggregv}
Let $\{\phi(\theta),  \theta\in T_0\}$ be an admissible family
of random
variables which is RCE. Suppose that $E [\operatorname{ess}\sup
_{\theta\in
T_0}\phi(\theta)]<\infty$.

There then exists an RCLL supermartingale $(v_t)$ which aggregates the
family $\{v(S), S\in T_0\}$ defined by (\ref{vs}),
that is, for each
stopping time $S$, $ v(S)=v_{S}$  \mbox{ a.s.}\vspace*{-1pt}
\end{Proposition}

\begin{pf}
The family
$\{v(S) , S\in T_0\}$ is a supermartingale system (Proposition~\ref{P1.SuperM}) and has the RCE property (Proposition \ref{L1}). The
result clearly follows by applying the aggregation property of
supermartingale systems (Proposition \ref{P.SMA}).\vspace*{-1pt}
\end{pf}

\begin{Theorem} \label{P.prop} Suppose the reward is given by an RC
and LCE admissible family
$\{ \phi(\theta), \theta\in T_0\}$ such that $E [\operatorname
{ess}\sup
_{\theta\in
T_0}\phi(\theta)]<\infty$.

Let $(\phi_t)$ be the progressive process given by Theorem 4.1 that
aggregates this family. Let $\{v(S), S\in T_0\}$ be the family of value
functions defined by~(\ref{vs}), and let $(v_t)$ be an RCLL adapted
process that aggregates the family $\{v(S), S\in T_0\}$.

The random variable defined by
%
\begin{equation}\label{thetaop}
\overline\theta(S)=\inf\{t\geq S,
v_t=\phi_t\}
\end{equation}
is the minimal optimal stopping time for $v(S)$, that is, $\overline
\theta (S) = \theta^* (S)$ a.s.\vspace*{-1pt}
\end{Theorem}

As for Theorem \ref{T.1}, the proof relies on the construction of a
family of stopping times that are approximatively optimal. The details,
which require some fine techniques of the general theory of processes,
are given in Section~\ref{3.14}.\vspace*{-1pt}

\begin{Remark}
In the case of an RCLL reward process supposed to be LCE, the above
theorem corresponds to the classical existence result [see El Karoui
(\citeyear{EK}) and Karatzas and Shreve (\citeyear{KS2})].\vspace*{-1pt}
\end{Remark}

\subsection{The optimal multiple stopping time problem}\label{3.13}
For simplicity, we study only the case when $d=2$. We will now prove
that the minimal optimal pair of stopping times $(\tau_1^*, \tau_2^*)$
defined by (\ref{mio}) can also be given in terms of \textit{hitting
times}. In order to do this, we first need to aggregate the value
function and the new reward.\vspace*{-1pt}

\subsubsection{Aggregation of the value function}\label{3.13.1}\vspace*{-1pt}

\begin{Proposition}\label{P2.Aggregv}
Suppose the reward is given by an RCE biadmissible family
$\{\psi(\theta,S), \theta,S\in T_0\}$ such that ${
E[\operatorname{ess}\sup_{\theta,S\in T_0}\psi(\theta,S)]<\infty
}$.\vadjust{\eject}

There then exists a supermartingale $( v_t)$ with RCLL paths that
aggregates the family $\{v(S), S\in T_0\}$ defined by (\ref{vS}), that
is, such for each $S$ $\in$ $T_0$, $v(S)= v_S $ a.s.
\end{Proposition}

\begin{pf}
The RCE property of $\{v(S),S\in T_0\}$ shown in Proposition \ref{L2.vRCE}, together with the supermartingale property [Proposition
\ref{P2.1}(3)] gives, by Proposition~\ref{P.SMA}, the desired result.
\end{pf}

\subsubsection{Aggregation of the new reward}\label{3.13.2}

We will now study the aggregation problem of the new reward family
$\{\phi(\theta),\theta\in T_0\}$.
Let us introduce the following definition.
\begin{Definition}
A biadmissible family $\{\psi(\theta,S),\theta
,S\in T_0\} $ is said to be \textit{uniformly right-continuous along
stopping times (URC)} if $ E[\operatorname{ess}\sup_{\theta, S\in T_0 } \psi (\theta,\break S) ] <
\infty$ and if for each nonincreasing sequence of stopping times
$(S_n)_{n \in\mathbb{N}}$ in $T_S$ which converges a.s.~to a stopping
time $S\in T_0$,
\[
\lim_{n\to\infty}\Bigl[\operatorname{ess}\sup_{\theta\in T_S}\{ |\psi(\theta,S_n) - \psi(\theta,S)| \}\Bigr]=0 \qquad\mbox{a.s.}
\]
and
\[
\lim_{n\to\infty}\Bigl[\operatorname{ess}\sup_{\theta\in T_S} \{|\psi(S_n, \theta) - \psi(S, \theta)| \}\Bigr]=0\qquad\mbox{a.s.}
\]
\end{Definition}

The following right continuity property holds true for the new reward
family.

\begin{Theorem} \label{T2.NGRC}
Suppose that the admissible family of positive random variables
$\{\psi(\theta,S),\theta,S\in T_0\}$ is URC. The family of positive
random variables $\{\phi(S), S\in T_0\}$ defined by (\ref{phi}) is then
RC.
\end{Theorem}

\begin{pf}
As $\phi(\theta)=\max[u_1(\theta),u_2(\theta
)]$, it is sufficient to show the RC property for the family $\{
u_1(\theta),\theta\in T_0\}$.

Now, for all $\theta\in T_0$, $u_1(\theta)=U_1(\theta,\theta)$ a.s.,
where
%
\begin{equation}\label{U1bis}
U_1 (\theta, S) = \operatorname{ess}\sup_{\tau_1 \in T_{\theta}}
E[\psi(\tau_1,
S)
|\mathcal{F}_{\theta}]\qquad\mbox{a.s.}
\end{equation}
Hence, it is sufficient to prove that $\{U_1 (\theta, S), \theta,
S \in T_0\}$ is RC.

Let $\theta,S\in T_0$ and $(\theta_n)_n$, $(S_n)_n$ be nonincreasing
sequences of
stopping times in $T_0$ that converge to $\theta$ and $S$ a.s. We have
\[
| U_1 (\theta, S)-U_1 (\theta_n, S_n) | \leq\underbrace{| U_1
(\theta, S) - U_1 (\theta_n , S) |}_{\mbox{(I)}} + \underbrace{| U_1
(\theta_n, S) - U_1 (\theta_n, S_n) |}_{\mbox{(II)}}.
\]

(I) \textit{tends to} 0 \textit{as} $n\to\infty$.

For each $S\in T_0$, as $\{\psi(\theta,S),\theta\in T_0\}$ is an
admissible
family of positive random variables which is RC, Proposition
\ref{P2.Aggregv} gives the existence of an RCLL adapted process
$(U^{1,S}_t)$ such that for each stopping time $\theta$ $\in$ $T_0$,
%
\begin{equation}\label{tuu}
U^{1,S}_{\theta}=U_1 (\theta, S) \qquad\mbox{a.s.}
\end{equation}
(I) can be rewritten as $
| U_1 (\theta, S) - U_1 (\theta_n, S) | =
| U^{1,S}_{\theta} - U^{1,S}_{\theta_n} |
$ a.s.,\vspace*{-2pt} which converges a.s.~to $0$ as $n$ tends to $\infty$ by the
right continuity of the process $(U^{1, \theta}_t)$.

(II) \textit{tends to} $0$ \textit{as} $n\to\infty$.

By definition of the value function $U_1(\cdot,\cdot)$ (\ref{U1bis}), it
follows that
\begin{eqnarray*}
| U_1 (\theta_n, S) - U_1 (\theta_n , S_n) |
& \leq&
E\Bigl(\operatorname{ess}\sup_{\tau_1 \in T_{\theta_n} } | \psi(\tau_1,S) - \psi(\tau_1, S_n) |\bigl|\mathcal{F}_{\theta_n}\Bigr )
\\
& \leq&
E( Z_m  |  \mathcal{F}_{\theta_n} )\qquad\mbox{a.s.}
\end{eqnarray*}
for any $n \geq m$, where $ Z_m := \sup_{r \geq m }
\{ \operatorname{ess}\sup_{\tau\in T_0 } | \psi(\tau, S_r)
- \psi(\tau, S)
| \}$ and $( E( Z_m  | \mathcal{F}_{t}))_{t
\geq
0}$ is an RCLL version of the conditional expectation. Hence, by the
right continuity of this process, for each fixed $m \in\mathbb{N}$, the
sequence of random variables $ ( E( Z_m  |  \mathcal{F}_{\theta_n}) )_{n \in\mathbb{N}}$ converges a.s.~to $ E(
Z_m  |
\mathcal{F}_{\theta})$ as $n$ tends to $\infty$. It follows that for each
$m \in \mathbb{N}$,
%
\begin{eqnarray}\label{UZ}
\limsup_{n \to\infty} | U_1 (\theta_n, S) - U_1 (\theta_n , S_n) | &
\leq& E( Z_m  |  \mathcal{F}_{\theta} )
\qquad
\mbox{a.s.}
\end{eqnarray}
Now, the sequence $(Z_m)_{m \in\mathbb{N}}$ converges a.s.~to $0$ and
\[
\vert Z_m \vert \leq 2 \operatorname{ess}\sup_{\theta,S\in
T_0}    \psi (\theta,S)   \qquad
\mbox{a.s.}
\]
Note that the second member of this inequality is integrable. By the
Lebesgue theorem for the conditional expectation, $ E( Z_m  |
\mathcal{F}_{\theta} )$ converges to $0$ in $L^1$ as~$m$ tends to
$\infty$. The sequence $(Z_m)_{m \in\mathbb{N}}$ is decreasing. It follows
that the sequence $\{ E( Z_m  |  \mathcal{F}_{\theta} )
\}_{m \in\mathbb{N}}$ is also decreasing and hence converges a.s.
Since this
sequence converges to $0$ in $L^1$, its limit is also $0$ almost
surely. By~let\-ting~$m$ tend to $\infty$ in (\ref{UZ}), we obtain
\[
\limsup_{n \to\infty} | U_1 (\theta_n, S) - U_1 (\theta_n ,S_n) |
\leq 0
\qquad
\mbox{a.s.}
\]
The proof of Theorem \ref{T2.NGRC} is thus complete.
\end{pf}

\begin{Corollary}[(Aggregation of the new reward)]\label{C2.1}
Under the same hypothesis as Theorem \ref{T2.NGRC}, there exists some
progressive right-continuous adap\-ted process $(\phi_t )$ which
aggregates the family $\{ \phi(\theta),$ $\theta\in T_0 \}$,
that is, $ \phi_{\theta}= \phi(\theta)$ a.s.~for each $\theta\in T_0$,
and such that there exists a decreasing sequence of right-continuous
processes $(\phi_t ^{n})_{n\in\mathbb{N}}$ that converges to
$(\phi_t)$.
\end{Corollary}

\begin{pf}
This follows from the right continuity of the new
reward (Theorem~\ref{T2.NGRC}) which we can aggregate (Theorem
\ref{T.Aggregphi}).
\end{pf}

\begin{Remark}
For the optimal $d$-stopping\vspace*{1pt} time problem, the same result holds for
URC $d$-admissible families $\{\psi(\theta)$, $\theta\in T_0^d\}$,
that is, families that satisfy $E[\operatorname{ess}\sup_{\theta\in T_0}   \psi (\theta)] < \infty$ and
\[
\lim_{n\to\infty} \operatorname{ess}\sup_{\theta\in T_0} \bigl|\psi
^{(i)}(\theta,S)-\psi
^{(i)}(\theta,S_n)\bigr|=0
\]
for $i=1,\ldots,d,\theta,S\in T_0$ and sequences $(S_n)$ in $T_0$
such that
$S_n\downarrow S$ a.s.

The proof is strictly the same, with $U_1(\theta,S)$ replaced by
$U^{(i)}(\theta,S)$ for~$\theta$, $S\in T_0$ and $\psi(\tau,S)$ with
$\tau,S\in T_0$
replaced by $\psi^{(i)}(\tau,S),$ with $\tau\in T_0^{d-1}$ and $S\in T_0$.
\end{Remark}

\subsubsection{Optimal multiple stopping times as hitting times of processes}\label{3.13.3}

As before, for the sake of simplicity, we suppose that $d=2$. Suppose
that $\{\psi(\theta,S), \theta, S\in T_0\}$ is a URC and ULCE
biadmissible family.
Let $\{\phi(\theta),\theta\in T_0\}$ be the new reward family. By Theorem
\ref{T2.NRcadgE}, this family is LCE. Furthermore, by Theorem
\ref{T2.NGRC}, this family is RC. Let $(\phi_t)$ be the progressive
process that aggregates this family, given by Theorem
\ref{T.Aggregphi}. Let $(u_t)$ be an RCLL process that aggregates the
value function associated with $(\phi_t)$. By Theorem \ref{P.prop}, the stopping time
\[
\theta^ *=\inf\{t\geq S, u_t=\phi_t\}
\]
is optimal for $u(S)$.

The family $\{\psi(\theta,\theta^*)$, $\theta\in T_{\theta^*}\}$ is
admissible, RC and LCE. Let $(\psi^1_t)$ be the progressive process
that aggregates this family given by Theorem \ref{T.Aggregphi}. Let
$(v^1_t)$ be an RCLL\vspace*{1pt} process that aggregates the value function
associated with $(\psi^1_t)$. By Theorem \ref{P.prop} the stopping time
$\theta^*_1=\inf\{t\geq\theta^*, v^1_t=\psi^1_t\}$ is optimal for
$v^1_{\theta^*}$ and $v^1_{\theta^*}=u^1(\theta^*)$.\vadjust{\goodbreak}

The family $\{\psi(\theta^*,\theta)$, $\theta\in T_{\theta^*}\}$ is
admissible, RC and LCE. Let $(\psi^2_t)$ be the progressive process
that aggregates this family given by Theorem \ref{T.Aggregphi}. Let
$(v^2_t)$ be an RCLL process\vspace*{1pt} that aggregates the value function
associated with $(\psi^2_t)$. By Theorem \ref{P.prop}, the stopping
time $\theta^*_2=\inf\{t\geq\theta^*, v^2_t=\psi^2_t\}$ is optimal
for $v^2_{\theta^*}$, and $v^2_{\theta^*}=u_2(\theta^*)$.

By Proposition \ref{Pconstruction}, the pair of stopping times
$(\tau_1^*,\tau_2^*)$ defined by
%
\begin{equation}\label{eto}
\tau_1^*=\theta^{*}{\mathbf{1} }_{B}+ \theta^{*}_1{\mathbf{1} }_{B^c},
\qquad
\tau_2^*=\theta^*_2{\mathbf{1} }_{B}+ \theta^{*}{\mathbf{1} }_{B^c},
\end{equation}
where $B=\{u_1(\theta^*) \leq u_2(\theta^ *)\}= \{v^1_{\theta^*} \leq v^2_{\theta^ *}\}$,
is optimal for $v(S)$.

\begin{Theorem}
Let $\{\psi(\theta,S), \theta,S\in T_0\}$ be a biadmissible
family which is
URC and ULCE.
The pair of stopping times $(\tau_1^*, \tau_2^*)$ defined by (\ref{eto}) is then optimal for $v(S)$.
\end{Theorem}

Note that the above construction of $(\tau_1^*,\tau_2^*)$ as hitting
times of processes requires stronger assumptions on the reward than
those made in Theorem~\ref{Topt}. Furthermore, let us emphasize that
it also requires some new aggregation results (Theorems
\ref{T.Aggregphi} and~\ref{P.prop}).

\subsection{\texorpdfstring{Proofs of Proposition \protect\ref{P.SMA} and
Theorems \protect\ref{T.Aggregphi} and \protect\ref{P.prop}}{Proofs of Proposition 4.1 and
Theorems 4.1 and 4.2}}
\label{3.14}

We now give the proofs of Proposition \ref{P.SMA} and Theorems
\ref{T.Aggregphi} and~\ref{P.prop}.

First, we give the short proof of the classical Proposition \ref{P.SMA}
which we recall here (for the reader's convenience).

\renewcommand{\theProposition}{4.1}
\begin{Proposition}
Let $\{h(S), S\in T_0\}$ be a supermartingale system which satisfies $h(0) < \infty$ and which is
right-continuous along stopping times in expectation. There then exists
an RCLL adapted process $(h_t)$ which
aggregates the family $\{h(S), S \in T_{0}\}$, that is, $h_S =
h(S)$ a.s.
\end{Proposition}

\begin{pf}
Let us consider the process $( h(t)
)_{0 \leq t \leq T}.$ It is a supermartingale and the function $t
\mapsto E(h(t))$ is right-continuous. By classical results [see Theorem~3.13 in Karatzas and Shreve (\citeyear{KS1})], there exists an RCLL
supermartingale $( h_t)_{0 \leq t \leq T}$ such that for each $t \in
[0,T]$, $h_t = h(t)$ a.s. It is then clear that for each dyadic
stopping time $S$ $\in$ $T_0$, $h_S = h(S)$ a.s. (for details, see Part
2 of the proof of Theorem \ref{T.1}). This implies that
%
\begin{equation}\label{eta}
E[h_S]= E[h(S)].
\end{equation}
Since the process $( h_t)_{0 \leq t \leq T}$ is RCLL and since the
family $\{ h(S), S \in T_0\}$ is right-continuous in expectation,
equality (\ref{eta}) still holds for any stopping time $S$ $\in$
$T_0$. It then remains to show that $h_S = h(S)$ a.s., but this is
classical. Let $A \in\mathcal{F}_S$ and define $S_A = S {\mathbf{1}}_A + T
{\mathbf{1}}_{A^c}$. Since $S_A$ is a stopping time, $E[h_{S_A}]=
E[h(S_A)]$. Since $h_T = h(T)$ a.s., it gives that $E[h_S {\mathbf{1}}_A ] =
E[h(S){\mathbf{1}}_A ],$ from which the desired result follows.
\end{pf}


We now give the proof of Theorem \ref{T.Aggregphi}.

\renewcommand{\theTheorem}{4.1}
\begin{Theorem}Suppose that
the admissible family of positive random variables $\{ \phi(\theta),
\theta\in T_0 \}$ is right-continuous along stopping times. There
then exists a progressive process $( \phi_t ) $ such that for each
$\theta$ $\in$ $T_0$, $ \phi_{\theta}= \phi(\theta)  \mbox{ a.s.} $
and such that there exists a nonincreasing sequence of
right-continuous processes~$(\phi_t ^{n})_{n\in\mathbb{N}}$ such
that for each
$(\omega, t)$ $\in$ $\Omega\times[0,T]$, $\lim_{n \rightarrow
\infty}
\phi_t ^{n}(\omega) = \phi_t (\omega)$.
\end{Theorem}


\begin{pf}
For each $n$ $\in$ $\mathbb{N}^*$, let us
define a
process $(\phi_t ^{n})_{t \geq0}$ that is a function of $(\omega,t)$
by
%
\begin{equation}\label{trois}
\phi^{n} _t (\omega)=\sup_{s \in\mathbb{D}\cap]t, ([2^n t]+1)/2^n [
} \phi(s\wedge T)
\end{equation}
for each $(\omega, t)$ $\in$ $\Omega\times[0,T],$ where $ \mathbb
{D}$ is the
set of dyadic rationals.

For each $t \in[0,T]$ and each $\varepsilon> \frac{1}{2^n}$, the
process $(\phi_t ^{n})$ is $(\mathcal{F}_{t+ \varepsilon})$-adapted and,
for each $\omega \in \Omega$, the function $t \mapsto\phi_t
^n(\omega)$ is right-continuous. Hence, the process $( \phi_t ^n)$ is
also $(\mathcal{F}_{t+ \varepsilon})$-progressive. Moreover, the sequence
$(\phi_t ^{n})_{n \in\mathbb{N}^*}$ is\vadjust{\eject} decreasing. Let $\phi_t$ be
its limit,
that is, for each $(\omega, t) \in \Omega\times[0,T]$,
\[
\phi_t (\omega)= \lim_{n \to\infty} \phi^{ n}_t (\omega).
\]
It follows that for each $\varepsilon> 0$, the process $(\phi_t )$ is
$(\mathcal{F}_{t+ \varepsilon})$-progressive. Thus, $(\phi_t )$ is
$(\mathcal{F}_{t^+})$-progressive and consequently $(\mathcal{F}_{t})$-progressive
since $\mathcal{F}_{t^+}= \mathcal{F}_{t}$.

\textit{Step 1}: Fix $\theta\in T_0$. Let us show that
$\phi_\theta\leq\phi(\theta)$ a.s.

Let us suppose, by contradiction, that the above inequality does not
hold. There then exists $\varepsilon>0$ such that the set $A=\{ \phi
(\theta)\leq\phi_\theta-\varepsilon\}$ satisfies \mbox{$P(A)>0$}.

Fix $n \in N$. For all $\omega\in A,$ we have that
$\phi(\theta)(\omega) \leq\phi^{n}_{\theta(\omega)}(\omega)
-\varepsilon,$ where $\phi^{n}_{\theta(\omega)}(\omega) $ is
defined by
(\ref{trois}) with $t$ replaced by $\theta(\omega)$.

By definition of $\phi^{n}$ there exists $t\in\,]\theta(\omega),
\frac{[2^n\theta(\omega)]+1}{2^n}[\,\cap\,\mathbb{D}$ such that
\[
\phi(\theta)(\omega) \leq\phi(t)(\omega) -\frac{\varepsilon}{2}.
\]
We introduce the following subset of $[0,T]\times\Omega$:
\[
\overline A_n=\biggl\{ (t,\omega)  , t\in\,\biggr]\theta(\omega), \frac
{[2^n\theta(\omega)]+1}{2^n}\biggl[\,\cap\,\mathbb{D}\mbox{ and } \phi
(\theta
)(\omega) \leq\phi(t)(\omega) -\frac{\varepsilon}{2} \biggr\}.
\]
First, note that $\overline A_n$ is optional. Indeed, we have
$\overline A_n= \bigcup_{t \in\mathbb{D}} \{t\} \times
B_{n,t},$
where
\[
B_{n,t}=\biggl\{ \theta< t < \frac{[2^n\theta]+1}{2^n}\biggr\} \cap
\biggl\{
\phi(\theta) \leq\phi(t) -\frac{\varepsilon}{2} \biggr\},
\]
and the
process $(\omega, t) \mapsto{\mathbf{1}}_{B_{n,t}} (\omega)$ is optional
since $\theta$ and $\frac{[2^n\theta]+1}{2^n}$ are stopping times and
$\{ \phi(\theta), \theta\in T_0 \}$ is admissible. Also, $A$ is
included in $\pi(\overline A_n)$, the projection of $\overline A_n$
onto $\Omega$, that is,
\[
A\subset\pi(\overline A_n) =\{ \omega\in\Omega,\ \exists t \in[0, T]  \mbox{ s.t. } (t, \omega) \in\overline A_n
 \}.
\]
Hence, by a section theorem [see Dellacherie and Meyer (\citeyear{DM1}),
Chapter IV], there exists a dyadic stopping time $T_n$ such that for
each $\omega$ in $\{T_n < \infty\}$, $(T_n(\omega), \omega)\in\overline A_n$ and
\[
P(T_n < \infty) \geq P(\pi(\overline A_n)) - \frac{ P(A)}{ 2^{n+1}}
\geq P(A) - \frac{ P(A)}{ 2^{n+1}}.
\]
Hence, for all $\omega$ in $\{T_n < \infty\}$
\[
\phi(\theta)(\omega) \leq\phi(T_n (\omega))-\frac{\varepsilon
}{2}\quad \mbox{and}\quad  T_n(\omega)\in\biggl]\theta(\omega), \frac
{[2^n\theta(\omega)]+1}{2^n}\biggr[\cap\mathbb{D}.
\]
Note that
\[
P\biggl( \bigcap_{n \geq1} \{T_n < \infty\}\biggr) \geq P(A) - \biggl( \sum_{n
\geq1} \frac{ P(A)}{ 2^{n+1}} \biggr)
\geq\frac{P(A)} {2} >0.
\]

Put $\overline T_n= T_1\wedge\cdots\wedge T_n$. We have $\overline
T_n \downarrow\theta$ and $\phi(\theta) \leq\phi(\overline
T_n)-\frac{\varepsilon}{2}$ for each $n$ on $\bigcap_{n \geq1} \{T_n <
\infty\} $. By letting $n$ tend to $\infty$ in this inequality, since
$\{ \phi(\theta), \theta\in T_0 \}$ is right-continuous along
stopping times, we derive that $\phi(\theta) \leq
\phi(\theta)-\frac{\varepsilon}{2}$ a.s.~on $\bigcap_{n \geq1} \{T_n <
\infty\}$, which gives the desired contradiction.

\textit{Step 2}: Fix $\theta$ $\in$ $T_0$. Let us show that
$\phi(\theta) \leq\phi_\theta $ a.s.

Put $T^n= \frac{[2^n\theta]+1}{2^n}$. The sequence $(T^n)$ is a
nonincreasing sequence of stopping times such that $ T^n\downarrow
\theta$. Moreover, note that since the family $\{ \phi(\theta),\break
\theta\in T_0 \}$ is admissible, for each $d$ $\in\mathbb{D}$,
for almost
every $\omega\in\{ T^{n+1}= d \}$, $\phi(T^{n+1})(\omega) =
\phi(d)(\omega)$. Now, we have $T^{n+1}\in\,]\theta, T^n[\,\cap\,\mathbb{D}$.
Also, for each \mbox{$\omega$ $\in$ $\Omega$} and each $d\in\,]\theta(\omega), T^n(\omega)[\,\cap\,\mathbb{D}$,
\[
\phi(d)(\omega) \leq\sup_{s\in\,]\theta(\omega), T^n(\omega)[\cap
\mathbb{D}} \phi(s)(\omega)= \phi^{n}_{\theta(\omega)} (\omega),
\]
where the last equality follows by the definition of
$\phi^{n}_{\theta(\omega)} (\omega)$ [see (\ref{trois}), with~$t$
replaced by $\theta(\omega)$]. Hence,
\[
\phi(T^{n+1}) \leq\phi^{n}_\theta\qquad\mbox{a.s.}
\]
Letting $n$ tend to $\infty$, by using the right-continuous property of
$\{ \phi(\theta), \theta\in T_0 \}$ along stopping times and the
convergence of $\phi^{n}_{\theta(\omega)} (\omega)$ to
$\phi_{\theta(\omega)} (\omega)$ for each~$\omega$, we derive that
$\phi(\theta)\leq\phi_\theta$ a.s.
\end{pf}

We now give the proof of Theorem \ref{P.prop}.

\renewcommand{\theTheorem}{4.2}
\begin{Theorem}
$\!\!\!\overline\theta (S)\!=\!\inf\{t\geq S,  v_t\!=\!\phi_t\} $ is an optimal
stopping time~for~$v_S$.
\end{Theorem}

\begin{pf}
We begin by constructing a family of stopping
times that are approximatively optimal. For $\lambda\in\,]0,1[$,
define the stopping time
%
\begin{equation} \label{otlS}
\overline\theta^{\lambda} (S) := \inf\{ t \geq S  ,  \lambda
v_t \leq\phi_t \} \wedge T.
\end{equation}
The proof follows the proof of Theorem \ref{T.1} exactly, except for
Step 1, which corresponds to the following lemma.

\begin{Lemma} \label{lemo}
For each $S\in T_0$ and $\lambda\in\,]0,1[$,
%
\begin{equation}\label{lemu}
\lambda v_{\overline\theta^{\lambda} (S) } \leq\phi_{\overline
\theta^{\lambda} (S) }\qquad \mbox{a.s.}
\end{equation}
\end{Lemma}

By the same arguments as in the proof of Theorem \ref{T.1}, $\overline
\theta^{\lambda} (S)$ is nondecreasing with respect to $\lambda$ and
converges as $\lambda\uparrow1$ to an optimal stopping time which
coincides with $\overline\theta(S)$ a.s.
\end{pf}

\begin{pf*}{Proof of Lemma \ref{lemo}}
To simplify notation,
$\overline\theta^{\lambda} (S)$ will be written as~$\overline\theta
^{\lambda}$. For the sake of simplicity, without loss of generality,
we suppose that $t \mapsto v_t(\omega)$ is RCLL for each $\omega\in
\Omega$.

Fix $\omega\in\Omega$. In the following, we use only simple analytic
arguments.

By definition of $\overline\theta^{\lambda} (\omega)$ (\ref{tlS}),
for each $n \in
\mathbb{N}^*$,
there exists $t$ $ \in[ \overline\theta^{\lambda}(\omega),
\overline\theta^{\lambda}(\omega) + \frac{1 } {n} [$ such that $
\lambda v_t (\omega) \leq\phi_t (\omega).$

Also, note that for each $m \in\mathbb{N}^*$, $ \phi_t (\omega)
\leq\phi
_t ^m (\omega)$.

Now, fix $m \in\mathbb{N}^*$ and $\alpha>0$.

By the right continuity of $t\mapsto v_t(\omega)$ and $t\mapsto
\phi^m_t(\omega)$, there exists $t^m_n(\omega)$ $\in$ $\mathbb
{D}\cap
[\overline\theta^{\lambda}(\omega),\overline\theta^{\lambda
}(\omega)
+ \frac{1 } {n} [$ such that
%
\begin{equation}\label{lemd}
\lambda v_{t^m_n(\omega)} (\omega) \leq\phi_{t^m_n(\omega)} ^m
(\omega)+ \alpha.
\end{equation}
Note that ${\lim_{n\to\infty}t^m_n(\omega)} =
\overline
\theta^{\lambda}(\omega)$ and $t^m_n(\omega) \geq\overline
\theta^{\lambda}(\omega)$ for any $n$.
Again, by using the right continuity of $t\mapsto
v_t(\omega)$ and $t\mapsto\phi^m_t(\omega),$ and by letting $n$ tend
to $\infty$ in (\ref{lemd}), we derive that
\[
\lambda v_{\overline\theta^{\lambda} (\omega) }(\omega) \leq
\phi^m_{\overline\theta^{\lambda} (\omega) }(\omega)+ \alpha,
\]
and this inequality holds for each $\alpha>0$, $m \in\mathbb{N}^*$ and
$\omega$ $\in$ $\Omega$. By letting $m$ tend to $\infty$ and
$\alpha$
tend to $0$, we derive that for each $\omega$ $\in$ $\Omega$,
$\lambda
v_{\overline\theta^{\lambda} (\omega)} (\omega) \leq\phi
_{\overline
\theta^{\lambda} (\omega) }(\omega),$ which completes the proof of the
lemma.
\end{pf*}

\begin{appendix}
\label{app}
\section{}\label{app1}
We recall the following classical theorem [see, e.g., Karatzas and
Shreve (\citeyear{KS2}), Neveu (\citeyear{Neveu})].

\renewcommand{\theTheorem}{A.1}
\begin{Theorem}[(Essential supremum)]\label{TA}
Let $(\Omega,\mathcal{F},P)$ be a probability space and let $
\mathcal{X}$ be a
nonempty family of positive random variables defined on
$(\Omega,\mathcal{F},P)$. There exists a random variable $X^*$ satisfying:
\begin{enumerate}
\item
for all $X\in\mathcal{X} $, $X\leq X^*$ a.s.;
\item
if $Y$ is a random variable satisfying $X\leq Y$ a.s.~for all $X\in
\mathcal{X} $, then $X^*\leq Y$~a.s.
\end{enumerate}
This random variable, which is unique a.s., is called the \textit{essential supremum of} $\mathcal{X}$ and is denoted $\operatorname
{ess}\sup
\mathcal{X}$.

Furthermore, if $\mathcal{X} $ is closed under pairwise maximization
\textup{(i.e.,} $X,Y\in\mathcal{X} $ implies $X\vee Y \in\mathcal{X}
$), then there is a nondecreasing sequence $\{Z_n\}_{n\in
\mathbb{N}
}$ of random variables in $\mathcal{X} $ satisfying $X^*=
{\lim_{n\to\infty} Z_n}$ a.s.
\end{Theorem}



\section{}\label{app2}

\subsection{Characterization of minimal optimal double stopping time}\label{app21}

In order to give a characterization of \textit{minimal optimal} stopping
times, we introduce the following partial order relation on $\mathbb{R}^2$:
$(a,b)\prec(a',b')$ if and only~if
\[
[(a\wedge b < a'\wedge b') \mbox{ or } (
a\wedge b = a'\wedge b'\mbox{ and } a\leq a'\mbox{ and } b\leq b')].
\]

Note that although the minimum of two elements of $\mathbb{R}^2$ is not
defined, the infimum, that is, the greatest minorant of the couple,
does exist and $\inf[ (a,b), (a',b')]= {\mathbf{1}}_{\{a\wedge b <
a'\wedge
b'\}}(a,b)+ {\mathbf{1}}_{\{a'\wedge b' < a\wedge b\}}(a',b') +
{\mathbf{1}}_{\{a\wedge b = a'\wedge b'\}}(a\wedge a' ,b\wedge b')$.

Note also that if $(\tau^*_1,\tau^*_2)$, $(\tau'_1,\tau'_2)$
$\in T_0 \times T_0$ are optimal for $v(S)$, then the infimum of the
couple $\inf[ (\tau^{*}_1,\tau^{*}_2), (\tau'_1,\tau'_2)]$, in
the sense of the relation $\prec$ a.s., is optimal for $v(S)$.

The two following assertions can be shown to be equivalent:
\begin{enumerate}[1.]
\item
a pair $(\tau^*_1,\tau^*_2)$ $\in T_0 \times T_0$ is \textit{minimal
optimal} for $v(S)$
(i.e, is the minimum for the order $\prec$ a.s.~of the set $\{
(\tau^*_1,\tau^*_2)\in T_S^2 $, $ v(S)= E[\psi
(\tau^*_1,\tau^*_2)|\mathcal{F}_S]\}$),
$\theta^*=\tau^*_1 \wedge\tau^*_2$ and $\theta^*_1$, $\theta^*_2$
$\in T_0$ are such that $\theta^*_2=\tau^*_2$ on $\{ \tau^*_1 < \tau
^*_2\}$ and
$\theta^*_1=\tau^*_1$ on $\{ \tau^*_1 > \tau^*_2\}$;
\item
\begin{enumerate}[(c)]
\item[(a)]$\theta^*$ $\in T_0$ is minimal optimal for $u(S)$;
\item[(b)]
$\theta^*_2$ $\in T_0$ is minimal optimal for $u_2({\theta^*})$ on
$\{u_1(\theta^*) < u_2(\theta^*)\};$
\item[(c)]
$\theta^*_1$ $\in T_0$ is minimal optimal for
$u_1({\theta^*})$ on $\{u_2(\theta^*) <  u_1(\theta^*)\}$,
and $\tau^*_1=\theta^*{\mathbf{1}}_{\{u_1(\theta^*)\leq
u_2(\theta^*)\}}+\theta^*_1{\mathbf{1}}_{\{u_1(\theta^*)> u_2(\theta
^*)\}}
$, $\tau^*_2=\theta^*{\mathbf{1}}_{\{u_2(\theta^*)\leq u_1(\theta^*)\}}+\theta^*_2\times\break{\mathbf{1}}_{\{u_2(\theta^*)> u_1(\theta^*)\}}$.
\end{enumerate}
\end{enumerate}

\subsection{Characterization of minimal optimal $d$-stopping times}\label{app22}

Consider the following partial order relation $\prec_d$ on $\mathbb{R}^d$
defined by induction in the following way:
for $d=1$, $\forall a,a'\in\mathbb{R}$, $a\prec_1 a'$
if and only if $a\leq a'$,
and
for $d>1$,
$
\forall(a_1,\ldots,a_d),(a_1',\ldots,a_d')\in\mathbb{R}^d,
(a_1,\ldots
,a_d)\prec_d (a_1',\ldots,a_d')$ if and only if
either $a_1\wedge\cdots\wedge a_d<a'_1\wedge\cdots
\wedge a'_d$ or
\[
\cases{a_1\wedge\cdots\wedge a_d=a'_1\wedge\cdots\wedge a'_d,\quad \mbox{and, for }i=1,\ldots, d,\cr
a_i=a_1\wedge\cdots\wedge a_d\quad\Longrightarrow\quad
\cases{a'_i=a'_1\wedge\cdots\wedge a'_d\quad \mbox{and}\cr
(a_1,\ldots,a_{i-1},a_{i+1},\ldots,a_d)\cr
\qquad\prec_{d-1}(a_1',\ldots,a'_{i-1},a'_{i+1},\ldots,a_d').
}
}
\]
Note that for $d=2$ the order relation $\prec_2$ is the order
relation $\prec$ defined above.

One can show that a $d$-stopping time $(\tau_1,\ldots,\tau_d)$ is the
$d$-minimal optimal stopping time for $v(S)$, that is, it is minimal
for the order $\prec_d$ in the set $\{ \tau\in T_S^d,  v(S)=E[
\psi(\tau)|\mathcal{F}_S] \}$ if and only if:
\begin{enumerate}
\item$\theta^*=\tau_1 \wedge\cdots\wedge\tau_d$ is minimal optimal
for $u(S)$;
\item
for $i=1,\ldots, d$,
$\theta^{*(i)}=\tau_i\in T_S^{d-1}$ is the $(d-1)$-minimal optimal
stopping time for $u^{(i)}({\theta^*})$ on the set
${\{u^{(i)}(\theta^*)\geq\bigvee_{k\neq i}
u^{(k)}(\theta^*)\}.}$
\end{enumerate}
\end{appendix}

\section*{Acknowledgments}
The authors thank Gilles Pag\`es and the anonymous referee for their
relevant remarks and suggestions.

%

\printaddresses


\begin{thebibliography}{13}

\bibitem[\protect\citeauthoryear{}{2008}]{CD}
%
\begin{barticle}[mr]
\bauthor{\bsnm{Carmona},~\bfnm{Ren{\'e}}\binits{R.}} \AND
\bauthor{\bsnm{Dayanik},~\bfnm{Savas}\binits{S.}}
(\byear{2008}).
\btitle{Optimal multiple stopping of linear diffusions}.
\bjournal{Math. Oper. Res.}
\bvolume{33}
\bpages{446--460}.
\bid{doi={10.1287/moor.1070.0301}, mr={2416002}}
\end{barticle}
%
\endbibitem

\bibitem[\protect\citeauthoryear{}{2008}]{CT}
%
\begin{barticle}[mr]
\bauthor{\bsnm{Carmona},~\bfnm{Ren{\'e}}\binits{R.}} \AND
\bauthor{\bsnm{Touzi},~\bfnm{Nizar}\binits{N.}}
(\byear{2008}).
\btitle{Optimal multiple stopping and valuation of swing options}.
\bjournal{Math. Finance}
\bvolume{18}
\bpages{239--268}.
\bid{doi={10.1111/j.1467-9965.2007.00331.x}, mr={2395575}}
\end{barticle}
%
\endbibitem

\bibitem[\protect\citeauthoryear{}{1975}]{DM1}
%
\begin{bbook}[vtex]
\bauthor{\bsnm{Dellacherie},~\bfnm{Claude}\binits{C.}} \AND
\bauthor{\bsnm{Meyer},~\bfnm{Paul-Andr{\'e}}\binits{P.-A.}}
(\byear{1975}).
\btitle{Probabilit\'es et Potentiel, Chap. I--IV},
\bedition{nouvelle} \'{e}dition.
\bpublisher{Hermann}, \baddress{Paris}.
\bid{mr={0488194}}
\bptnote{check year}
\end{bbook}
%
\endbibitem


\bibitem[\protect\citeauthoryear{}{1981}]{EK}
%
\begin{bincollection}[vtex]
\bauthor{\bsnm{El~Karoui},~\bfnm{N.}\binits{N.}}
(\byear{1981}).
\btitle{Les aspects probabilistes du contr\^ole stochastique}.
In
\bbooktitle{\'Ecole d'\'et\'e de Probabilit\'es de Saint-Flour IX-1979}.
\bseries{Lect. Notes in Math.}
\bvolume{876}
\bpages{73--238}.
\bpublisher{Springer}, \baddress{Berlin}.
\bid{mr={0637469}}%
\end{bincollection}
%
\endbibitem

\bibitem[\protect\citeauthoryear{}{1994}]{KS1}
%
\begin{bbook}[vtex]
\bauthor{\bsnm{Karatzas},~\bfnm{Ioannis}\binits{I.}} \AND
\bauthor{\bsnm{Shreve},~\bfnm{Steven~E.}\binits{S.~E.}}
(\byear{1994}).
\btitle{Brownian Motion and Stochastic Calculus},
\bedition{2nd} ed.
\bseries{Graduate Texts in Mathematics}
\bvolume{113}.
\bpublisher{Springer}, \baddress{New York}.
\bid{mr={1121940}}
\bptnote{check year}
\end{bbook}
%
\endbibitem

\bibitem[\protect\citeauthoryear{}{1998}]{KS2}
%
\begin{bbook}[mr]
\bauthor{\bsnm{Karatzas},~\bfnm{Ioannis}\binits{I.}} \AND
\bauthor{\bsnm{Shreve},~\bfnm{Steven~E.}\binits{S.~E.}}
(\byear{1998}).
\btitle{Methods of Mathematical Finance}.
\bseries{Applications of Mathematics (New York)}
\bvolume{39}.
\bpublisher{Springer}, \baddress{New York}.
\bid{mr={1640352}}
\end{bbook}
%
\endbibitem

\bibitem[\protect\citeauthoryear{}{2010}]{KQ}
%
\begin{bmisc}[vtex]
\bauthor{\bsnm{Kobylanski},~\bfnm{M.}\binits{M.}} \AND
\bauthor{\bsnm{Quenez},~\bfnm{M.~C.}\binits{M.~C.}}
(\byear{2010}).
Optimal multiple stopping in the Markovian case and
applications to finance. Working paper.
\end{bmisc}
%
\endbibitem

\bibitem[\protect\citeauthoryear{}{2010}]{KQR}
%
\begin{barticle}[mr]
\bauthor{\bsnm{Kobylanski},~\bfnm{Magdalena}\binits{M.}},
\bauthor{\bsnm{Quenez},~\bfnm{Marie-Claire}\binits{M.-C.}} \AND
\bauthor{\bsnm{Rouy-Mironescu},~\bfnm{Elisabeth}\binits{E.}}
(\byear{2010}).
\btitle{Optimal double stopping time problem}.
\bjournal{C. R. Math. Acad. Sci. Paris}
\bvolume{348}
\bpages{65--69}.
\bid{doi={10.1016/j.crma.2009.11.020}, mr={2586746}}
\end{barticle}
%
\endbibitem

\bibitem[\protect\citeauthoryear{}{1998}]{KR}
%
\begin{bmisc}[vtex]
\bauthor{\bsnm{Kobylanski},~\bfnm{M.}\binits{M.}} \AND
\bauthor{\bsnm{Rouy},~\bfnm{E.}\binits{E.}}
(\byear{1998}).
\btitle{Large deviations estimates for diffusion processes with Lipschitz
reflections}.
Th\`ese de Doctorat de L'universit\'e de Tours de M. Kobylanski
\bpages{17--62}.
\end{bmisc}
%
\endbibitem

\bibitem[\protect\citeauthoryear{}{1978}]{Ma}
%
\begin{bincollection}[mr]
\bauthor{\bsnm{Maingueneau},~\bfnm{M.~A.}\binits{M.~A.}}
(\byear{1978}).
\btitle{Temps d'arr\^et optimaux et th\'eorie g\'en\'erale}.
In \bbooktitle{S\'eminaire de {P}robabilit\'es, {XII} ({U}niv. {S}trasbourg,
{S}trasbourg, 1976/1977)}.
\bseries{Lecture Notes in Math.}
\bvolume{649}
\bpages{457--467}.
\bpublisher{Springer}, \baddress{Berlin}.
\bid{mr={0520020}}
\end{bincollection}
%
\endbibitem

\bibitem[\protect\citeauthoryear{}{1975}]{Neveu}
%
\begin{bbook}[vtex]
\bauthor{\bsnm{Neveu},~\bfnm{J.}\binits{J.}}
(\byear{1975}).
\btitle{Discrete-Parameter Martingales}, \bedition{revised} ed. \bseries{North-Holland Mathematical
Library}
\bvolume{10}.
\bpublisher{North-Holland}, \baddress{Amsterdam}.
\bnote{Translated from the French by T. P. Speed.}
\bid{mr={0402915}}%
\end{bbook}
%
\endbibitem

\bibitem[\protect\citeauthoryear{}{2006}]{Peskir}
%
\begin{bbook}[vtex]
\bauthor{\bsnm{Peskir},~\bfnm{Goran}\binits{G.}} \AND
\bauthor{\bsnm{Shiryaev},~\bfnm{Albert}\binits{A.}}
(\byear{2006}).
\btitle{Optimal Stopping and Free-Boundary Problems}.
\bpublisher{Birkh\"auser}, \baddress{Basel}.
\bid{mr={2256030}}
\end{bbook}
%
\endbibitem

\end{thebibliography}
\end{document}